\documentclass[11pt,letterpaper,reqno]{amsart}
\usepackage{amsmath,amsthm,amsfonts,amssymb,mathtools,algorithm,algpseudocode}
\usepackage{comment,bookmark,bm,color,tikz}
\usepackage{bm,makecell,multicol,graphicx}

\hypersetup{pdfstartview={FitH}}
\def\restrict#1{\raise-.5ex\hbox{\ensuremath|}_{#1}}
\newtheorem{Theorem}{{\bf Theorem}}[section]

\newtheorem{Corollary}[Theorem]{{\bf Corollary}}

\newtheorem{Remark}[Theorem]{{\bf Remark}}

\numberwithin{equation}{section}

\newcommand{\R}{\mathbb{R}}

\newcommand{\calS}{\mathcal{S}}
\newcommand{\calT}{\mathcal{T}}
\newcommand{\calX}{\mathcal{X}}

\newcommand{\range}{\text{range}}
\newcommand{\rank}{\text{rank}}

\begin{document}

\title[Randomized algorithms for coupled decompositions]{Randomized algorithms for coupled decompositions}
\author{Erna Begovi\'{c}~Kova\v{c}}\thanks{\textsc{Erna Begovi\'{c} Kova\v{c}}, University of Zagreb Faculty of Chemical Engineering and Technology, Maruli\'{c}ev trg 19, 10000 Zagreb, Croatia. \texttt{ebegovic@fkit.unizg.hr}}
\author{Anita Carevi\'{c}}\thanks{\textsc{Anita Carevi\'{c}}, University of Split Faculty of Electrical Engineering, Mechanical Engineering and Naval Architecture, Ru\dj era Bo\v{s}kovi\'{c}a 32, 21000 Split, Croatia. \texttt{carevica@fesb.hr}}
\author{Ivana \v{S}ain Glibi\'{c}}\thanks{\textsc{Ivana \v{S}ain Glibi\'{c}}, University of Zagreb Faculty of Science, Department of Mathematics, Bijeni\v{c}ka 30, 10000 Zagreb, Croatia. \texttt{ivanasai@math.hr}}

\thanks{This work has been supported by Croatian Science Foundation under the project UIP-2019-04-5200.}
\date{\today}

\renewcommand{\subjclassname}{\textup{2020} Mathematics Subject Classification}

\begin{abstract}
Coupled decompositions are a widely used tool for data fusion. As the volume of data increases, so does the dimensionality of matrices and tensors, highlighting the need for more efficient coupled decomposition algorithms. This paper studies the problem of coupled matrix factorization (CMF), where two matrices represented in low-rank form share a common factor. Additionally, it explores coupled matrix and tensor factorization (CMTF), where a matrix and a tensor are represented in low-rank form, also sharing a common factor matrix. We show that these problems can be solved using a direct approach with singular value decomposition (SVD), rather than relying on an iterative method. Knowing that matrices coming from real-world applications are often very large, the computational cost can be substantial. To address this issue and improve the efficiency, we propose new techniques for randomizing these algorithms. This includes a novel strategy for selecting a projection subspace that takes into account the contribution from both matrices involved in the decomposition equally.
We present extensive results of numerical tests that confirm the efficiency of our algorithms. Furthermore, as a novel approach and with a high success rate, we apply our randomized algorithms to the face recognition problem.
\end{abstract}

\maketitle

\section{Introduction}

Coupled decompositions of multiple data sets are broadly used in different engineering disciplines as a tool for data fusion. They are utilized in the analysis of data coming from different sources, for example, to better describe data obtained by different technologies or methods.
To name a few applications, coupled decompositions appear in chemometrics~\cite{Smilde_chem00,AcDuKoMo11,AcDuKo_chem11,EeDeLChem20}, signal processing~\cite{Smilde10,SoDeL13proc}, bioinformatics~\cite{Bagherian_bioinf20,Borsoi_bioinf23}, metabolimics~\cite{Acar_metabolimics15,Acar_metabolimics24a,Acar_metabolimics24b}, chromatography~\cite{Armstrong_chrom23}, etc.
In this paper, we apply them to the problem of face recognition. To the best of our knowledge, this is a novel approach. The standard procedure for the face recognition algorithm is to calculate the mean face image, derive the covariance matrix, extract its principal components, and use them to project the images onto a lower-dimensional subspace, as explained in~\cite{yang2004two,jia2017color,xiao2018two}, and generalized to tensors in~\cite{vasilescu2002multilinear,hached2021multidimensional,el2024tensor}. However, the coupled decomposition bypasses this procedure and achieves dimensionality reduction by extracting the common part of the images.

Coupled matrix factorization (CMF)~\cite{SiGoCMF08} decomposes a set of matrices in a way that they are represented in a low-rank format sharing one common factor. For two matrices $X\in\R^{m\times n_1}$ and $Y\in\R^{m\times n_2}$, their coupled rank-$k$ approximation is given by
\begin{equation}\label{CMF}
X\approx UV^T \quad \text{and} \quad Y\approx UW^T,
\end{equation}
where $U\in\R^{m\times k}$, $V\in\R^{n_1\times k}$, $W\in\R^{n_2\times k}$.

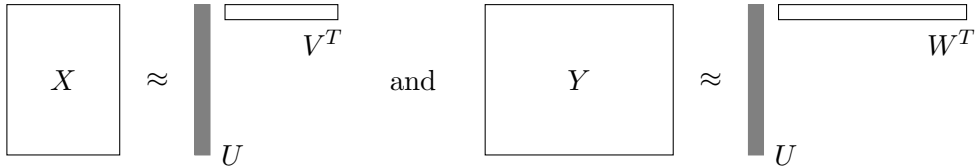
\begin{figure}[t]
\begin{tikzpicture}
\draw [draw=black] (0,0) rectangle (1.5,2);
\node at (2,1) {$\approx$};
\node at (0.75,1) {$X$};
\filldraw[fill=gray, draw=gray] (2.5,0) rectangle (2.7,2);
\draw[draw=black] (2.9,1.8) rectangle (4.4,2);
\node at (5.4,1) {\text{and}};
\node at (3,0) {$U$};
\node at (4.2,1.5) {$V^T$};
\end{tikzpicture}
\quad
\begin{tikzpicture}
\draw [draw=black] (0,0) rectangle (2.5,2);
\node at (3,1) {$\approx$};
\node at (1.25,1) {$Y$};
\filldraw[fill=gray, draw=gray] (3.5,0) rectangle (3.7,2);
\draw[draw=black] (3.9,1.8) rectangle (6.4,2);
\node at (4,0) {$U$};
\node at (6.2,1.5) {$W^T$};
\end{tikzpicture}
\caption{Graphical description of CMF.}
\label{fig:CMF}
\end{figure}

Big and complex data sets are often represented by tensors rather than matrices. Therefore, in addition to the matrix-matrix case, it is also important to consider tensor-matrix factorization, generalizing the problem~\eqref{CMF}. There are various algorithms for the coupled matrix and tensor factorization (CMTF). An often-used approach is the alternating least squares (ALS) method~\cite{Smilde_chem00,AcDuKo_chem11,AcDuKoMo11,SoDeL13proc,SoDeL15-2,Schenker21,Schenker21proc}, which is an iterative method that may encounter convergence issues. We are going to show how the CMTF of the tensor $\calX\in\R^{m\times n_2\times n_3}$ and matrix $Y\in\R^{m\times n}$ can be expressed in terms of CMF of the tensor matricization $X_{(1)}\in\R^{m\times n_2n_3}$ and matrix $Y$. Consequently, CMTF is simplified and solved using the singular value decomposition of the matrix $\begin{bmatrix} X_{(1)} & Y \end{bmatrix}\in\R^{m\times(n_2n_3+n)}$, instead of the alternating least squares technique.
This way, we do not obtain the full decomposition, as is the case with ALS, but we achieve the shared subspace faster and avoid potential convergence problems.
Since the aforementioned factorizations can be time-consuming for large-scale matrices and tensors, we look into a way to overcome this issue.

In recent years, randomized algorithms experienced a breakthrough in numerical linear algebra~\cite{HaMaTr11,Woo14,DriMa16,Ma16,MaTr20,MurrayRandNLA23}. Randomized algorithms are known to be significantly faster than deterministic algorithms, and they are reliable in many applications. In the context of coupled decompositions, randomized projections as an option for dealing with large data sets were briefly mentioned in~\cite{SoDeL15-2}. However, this paper thoroughly analyzes the potential of randomization in such decompositions, emphasizing the requirement of considering all matrices equally. We develop a randomized algorithm for coupled decomposition, inspired by the randomized SVD~\cite{HaMaTr11}, which has proven effective for a single matrix. Furthermore, we study more refined forms of randomization, similar to randomized subspace iteration~\cite{HaMaTr11,TrWe23} and randomized block Krylov iteration~\cite{TrWe23}, for the case of coupled decompositions. These approaches are tested on different examples, both for CMF and for CMTF, and the results are compared mutually, as well as with the non-randomized algorithms. In addition to the synthetic numerical examples, the algorithms are also tested on the face recognition problem. Given a database of faces, our randomized algorithms for coupled decompositions match a new face with a person from the database. Although this is not a traditional approach to face recognition, our algorithms demonstrated a very good
success rate.

In summary, the key contributions of this paper are as follows:
\begin{itemize}
\item A direct approach based on SVD to solve the coupled rank-$k$ approximation problem for two matrices and for matrix-tensor pairs is proposed.
\item A new strategy for constructing randomized projections in problems involving two matrices is developed. Based on this strategy, randomized algorithms for CMF and CMTF are created.
\item As a novel approach, our randomized algorithms for coupled decompositions are applied to the face recognition problem.
\end{itemize}

The paper is divided as follows. In Section~\ref{sec:CMF} we study coupled matrix factorization and introduce our basic algorithm (Algorithm~\ref{agm:CMF}). This algorithm is randomized in Subsection~\ref{ssec:randCMF}, while the numerical examples are given in Section~\ref{sec:CMF_num}. Then, we analyse coupled matrix and tensor factorization in Section~\ref{sec:CMTF}. The Tucker tensor decomposition is examined in Subsection~\ref{ssec:CMTF_Tucker} and CP tensor decomposition in Subsection~\ref{ssec:CMTF_CP}. The corresponding numerical examples are presented in Section~\ref{sec:CMTF_num}. In Section~\ref{sec:facerec} we apply our algorithms to the problem of face recognition. We end the paper with a short conclusion in Section~\ref{sec:conclusion}.

\section{Coupled matrix factorization}\label{sec:CMF}

We start with the coupled matrix factorization (CMF). Let $X\in\R^{m\times n_1}$ and $Y\in\R^{m\times n_2}$. The goal of CMF is to find a coupled rank-$k$ approximation of $X$ and $Y$ in the form given in the relation~\eqref{CMF}.
To solve this problem, we need to minimize the objective function
\begin{equation}\label{minM}
f(U,V,W)=\|X-UV^T\|_F^2+\|Y-UW^T\|_F^2\rightarrow\min.
\end{equation}

Before we construct an algorithm for solving the minimization problem~\eqref{minM}, we first show that the approximation problem~\eqref{CMF} is equivalent to the low-rank approximation of one matrix.

\begin{Theorem}\label{tm:cmf}
Let $X\in\R^{m\times n_1}$ and $Y\in\R^{m\times n_2}$.
Let $\Sigma_k\in\R^{k\times k}$ be a diagonal matrix of the largest $k$ singular values of the matrix $
\begin{bmatrix}
    X & Y
\end{bmatrix}\in\R^{m\times(n_1+n_2)}$, and let $U_k\in\R^{m\times k}$ and $V_k\in\R^{(n_1+n_2)\times k}$ contain left and right singular vectors of $\begin{bmatrix}
    X & Y
\end{bmatrix}$, respectively, corresponding to those singular values.
Then, solution of the minimization problem~\eqref{minM} is defined by $U_k$, $V_k$ and $\Sigma_k$.

Precisely, the objective function $f$ defined in~\eqref{minM} attains its minimum for $U=U_k$, $V$ equal to the first $n_1$ rows of $V_k\Sigma_k$, and $W$ equal to the remaining $n_2$ rows of $V_k\Sigma_k$.
\end{Theorem}

\begin{proof}
Let $U\in\R^{m\times k}$, $V\in\R^{n_1\times k}$, $W\in\R^{n_2\times k}$ be arbitrary matrices. For the fixed $X\in\R^{m\times n_1}$, $Y\in\R^{m\times n_2}$, using the properties of the Frobenius norm, we get
$$\|X-UV^T\|_F^2+\|Y-UW^T\|_F^2=\|\begin{bmatrix}
    X & Y
\end{bmatrix} - \begin{bmatrix}
    UV^T & UW^T
\end{bmatrix}\|_F^2=\|\begin{bmatrix}
    X & Y
\end{bmatrix} - UZ^T\|_F^2,$$
where $Z=\left[
  \begin{array}{c}
    V \\
    W
  \end{array}
\right]$.
Hence,
\begin{equation}\label{tm:eqcmf}
\min_{U,V,W}\left\{\|X-UV^T\|_F^2+\|Y-UW^T\|_F^2\right\}=\min_{U,Z}\|\begin{bmatrix}
    X & Y
\end{bmatrix} - UZ^T\|_F^2.
\end{equation}
The left-hand side in~\eqref{tm:eqcmf} is the best coupled rank-$k$ approximation of $X$ and $Y$, while the right-hand side corresponds to the best rank-$k$ approximation of the matrix $\begin{bmatrix}
    X & Y
\end{bmatrix}$.

It is well-known that a solution of the minimization problem on the right-hand side of~\eqref{tm:eqcmf} is given by the truncated SVD,
$$\begin{bmatrix}
    X & Y
\end{bmatrix}\approx U_k\Sigma_k V_k^T,$$
where $U_k, V_k,\Sigma_k$ are as in the statement of the theorem.
Then,
$$\min_{U,Z}\|\begin{bmatrix}
    X & Y
\end{bmatrix} - UZ^T\|_F^2 = \|\begin{bmatrix}
    X & Y
\end{bmatrix} - \bar{U}\bar{Z}^T\|_F^2,$$
where $\bar{U}=U_k$ and $\bar{Z}^T=\Sigma_kV_k^T$, that is, $\bar{Z}=V_k\Sigma_k$.

Therefore, it follows from~\eqref{tm:eqcmf} that
\begin{align*}
\min_{U,V,W}\left\{\|X-UV^T\|_F^2+\|Y-UW^T\|_F^2\right\} & =\|\begin{bmatrix}
    X & Y
\end{bmatrix} - \bar{U}\bar{Z}^T\|_F^2 \\
& = \|X-\bar{U}\bar{V}^T\|_F^2+\|Y-\bar{U}\bar{W}^T\|_F^2,
\end{align*}
for matrices $\bar{V}$ and $\bar{W}$ obtained by splitting the matrix $\bar{Z}$ into two parts, such that $\bar{V}$ contains the first $n_1$ rows of $Z$, while $\bar{W}$ contains the remaining $n_2$ rows of $\bar{Z}$.
\end{proof}

\begin{Corollary}
The best coupled rank-$k$ approximation of two matrices $X\in\R^{m\times n_1}$ and $Y\in\R^{m\times n_2}$ is equivalent to the best rank-$k$ approximation of the matrix $\begin{bmatrix}
    X & Y
\end{bmatrix}\in\R^{m\times(n_1+n_2)}$.
\end{Corollary}

\begin{proof}
The proof follows directly from Theorem~\ref{tm:cmf} and relation~\eqref{tm:eqcmf}.
\end{proof}

Now, using Theorem~\ref{tm:cmf} we can write the algorithm for solving the minimization problem~\eqref{minM}.
Algorithm~\ref{agm:CMF} is our baseline for the coupled matrix factorization. In the following subsection, we are going to incorporate different randomization techniques.

\begin{algorithm}
\caption{CMF --- Basic algorithm}\label{agm:CMF}
\renewcommand{\algorithmicrequire}{\textbf{Input:}}
\renewcommand{\algorithmicensure}{\textbf{Output:}}
\begin{algorithmic}
\Require $X\in\R^{m\times n_1}$, $Y\in\R^{m\times n_2}$, $k<\min\{n_1,n_2\}$
\Ensure $U\in\R^{m\times k}$, $V\in\R^{n_1\times k}$, $W\in\R^{n_2\times k}$
\State $[U_{XY},\Sigma,V_{XY}]=\texttt{svd}([X \ Y])$
\State $U=U_{XY}(:,1:k)$
\State $Z=V_{XY}(:,1:k)\Sigma(1:k,1:k)$
\State $V=Z(1:n_1,:)$
\State $W=Z(n_1+1:n_1+n_2,:)$
\end{algorithmic}
\end{algorithm}

\subsection{Randomized CMF}\label{ssec:randCMF}

The core of an efficient randomization technique for CMF is determining the projections of matrices $X$ and $Y$ onto the same subspace, as the accuracy of the result depends on it, i.e., how well it approximates the joint subspace of the direct sum
$\range(X) + \range(Y)$. Given Theorem \ref{tm:cmf}, a straightforward approach for selecting the projection matrix $Q$ is to compute the basis for the subspace $\range\left( \begin{bmatrix}
    X & Y
\end{bmatrix} \right)$.
However, we propose a more refined method. Specifically, we first determine the bases $Q_1$ and $Q_2$ for the subspaces $\range(X)$ and $\range(Y)$, respectively. Then, $Q$ is obtained by reorthogonalizing the columns of the matrix $\begin{bmatrix}
    Q_1 & Q_2
\end{bmatrix}$.
The approach involves generating projections of
$X$ and $Y$ onto the same subspace using a random Gaussian matrix $\Omega$:
\begin{equation}\label{randproj}
\widehat{X}=\Pi_{X\Omega}X, \quad \widehat{Y}=\Pi_{Y\Omega}Y.
\end{equation}
Subsequently, the coupled matrix factorization (CMF) of the reduced matrices
$\widehat{X}$ and $\widehat{Y}$, which are significantly smaller than
 $X$ and $Y$, is computed. The detailed procedure is outlined below.

For $X\in\R^{m\times n_1}$ and $Y\in\R^{m\times n_2}$ we generate random Gaussian matrices $\Omega_1\in{\R^{n_1\times k}}$ and $\Omega_2\in{\R^{n_2\times k}}$. Then we find the thin QR decompositions
\begin{equation}\label{Omega_rand}
Q_1R_1=X\Omega_1, \quad Q_2R_2=Y\Omega_2,
\end{equation}
such that $Q_1,Q_2\in\R^{m\times k}$. The matrices $Q_1$ and $Q_2$ form orthogonal bases of subspaces of $\range(X)$ and $\range(Y)$, respectively. To obtain a joint subspace of the direct sum $\range(X)+\range(Y)$, we reorthogonalize the columns of the matrix $\begin{bmatrix}
    Q_1 & Q_2
\end{bmatrix}$. This way we get a matrix $Q$ whose columns form an orthogonal basis of a subspace of $\range(X)+\range(Y)$. Matrix $Q$ has $m$ rows and $k+p$ columns, where $0\leq p\leq k$. Notice that $p<k$ if the subspaces of $\range(X)$ and $\range(Y)$ intersect. We can think of $p$ as an oversampling parameter.

When we have $Q$, we apply Algorithm~\ref{agm:CMF} on $\widehat{X}=Q^TX\in\R^{(k+p)\times n_1}$ and $\widehat{Y}=Q^TY\in\R^{(k+p)\times n_2}$. It returns $\widehat{U}\in\R^{(k+p)\times k}$, $V=\widehat{V}\in\R^{n_1\times k}$, $W=\widehat{W}\in\R^{n_2\times k}$. Finally, we set $U=Q\widehat{U}$.
Randomized CMF is presented in Algorithm~\ref{agm:randCMF}. In our code, we truncate $Q$ by reducing the number of columns from $2k$ to $k+p$, if possible, by examining the diagonal elements of the triangular matrix from the QR decomposition.

\begin{algorithm}
\caption{Randomized CMF}\label{agm:randCMF}
\renewcommand{\algorithmicrequire}{\textbf{Input:}}
\renewcommand{\algorithmicensure}{\textbf{Output:}}
\begin{algorithmic}
\Require $X\in\R^{m\times n_1}$, $Y\in\R^{m\times n_2}$, $k<\min\{n_1,n_2\}$
\Ensure $U\in\R^{m\times k}$, $V\in\R^{n_1\times k}$, $W\in\R^{n_2\times k}$
\State Generate random matrices $\Omega_1\in{\R^{n_1\times k}}$ and $\Omega_2\in{\R^{n_2\times k}}$.
\State $[Q_1,\sim]=\texttt{qr}(X\Omega_1,0)$
\State $[Q_2,\sim]=\texttt{qr}(Y\Omega_2,0)$
\State $[Q,\sim]=\texttt{qr}([Q_1,Q_2],0)$
\State $[\widehat{U},V,W]=\texttt{CMF}(Q^TX,Q^TY,k)$ \Comment{Algorithm~\ref{agm:CMF}}
\State $U=Q\widehat{U}$
\end{algorithmic}
\end{algorithm}

As can be seen from the above discussion and the relations~\eqref{randproj} and~\eqref{Omega_rand}, our randomization technique differs from the randomization of other matrix pair factorizations such as generalized singular value decomposition~\cite{wei2016tikhonov}. Since we are finding the joint subspace of a direct sum, our procedure includes both matrices to obtain the projection matrix $Q$.

We compare two approaches for generating projection matrices and compute
\begin{equation*}
    \|X-UV^T\|_F^2 + \|Y-UW^T\|_F^2,
\end{equation*}
which is the obtained minimum of the objective function~\eqref{minM}. The first approach is as in Algorithm~\ref{agm:randCMF}. In the second approach, the projection matrix $\bar{Q}$ is the orthogonal factor in the thin QR decomposition
\begin{equation}\label{Omega12}
\bar{Q}R=\begin{bmatrix}
    X & Y
\end{bmatrix}\bar{\Omega},
\end{equation}
where $\bar{\Omega}$ is $(n_1+n_2)\times k$ random Gaussian matrix, and CMF is performed on $\bar{Q}^TX$ and $\bar{Q}^TY$.

\begin{figure}[h!]
    \begin{minipage}{.5\textwidth}
    \centering
    \includegraphics[width=1\linewidth]{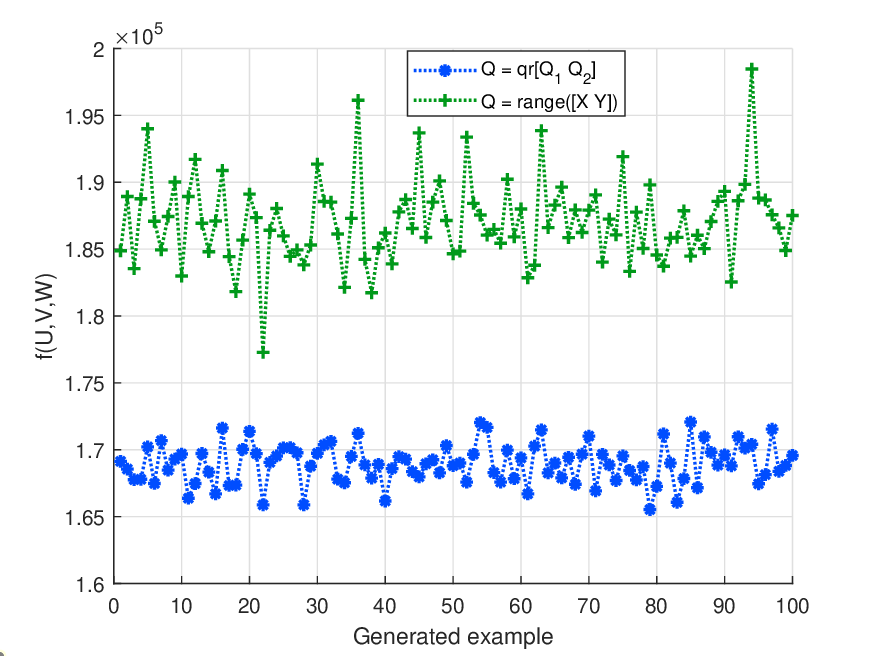}
    \end{minipage}%
    \begin{minipage}{.5\textwidth}
        \centering
        \includegraphics[width=1\textwidth]{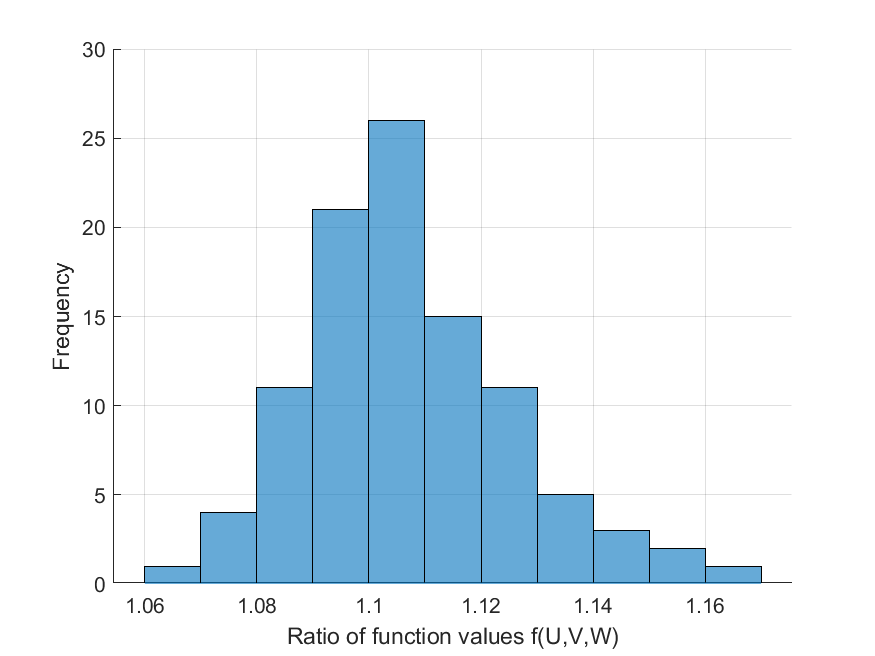}	
    \end{minipage}\\
    \caption{Comparison of two different approaches for computing the projection matrix $Q$ on $100$ randomly generated examples.}
    \label{fig:FirstExample}
\end{figure}

The following set of tests demonstrates that our approach is superior and can significantly impact the outcome of the randomization algorithms. We generate low-rank matrices $X\in\R^{500\times 200}$ and $ Y\in\R^{500\times 300}$
without imposing any specific structure on their singular values:
\begin{itemize}
     \item \texttt{X = rand($500$,$100$)*rand($100$,$200$);}
        \item \texttt{Y = rand($500$,$150$)*rand($150$,$300$);}
    \end{itemize}
The target rank for the coupled decomposition is fixed at $k=30$. We perform $100$ experiments on different, randomly generated, pairs of matrices $X$ and $Y$. The results are presented in Figure \ref{fig:FirstExample}. It is clear that, for every generated pair of matrices, our approach for choosing the projection matrix performs significantly better. In the majority of experiments, when the projection strategy is not applied, the value of the objective function $f(U,V,W)$ is around $10\%$ or more larger than in the case when our projection strategy is used.

We are going to also explain this behavior from the theoretical point of view. Firstly, note that the approximation error obtained by randomization can be estimated using the result for the randomized SVD from~\cite{HaMaTr11}.

\begin{Theorem}\label{tm:error}~\cite[Theorem 10.5]{HaMaTr11}
Suppose that $A$ is a real $m\times n$ matrix with singular values $\sigma_1\geq\sigma_2\geq\cdots\geq\sigma_{\min\{m,n\}}$. Choose a target rank $k\geq2$ and an oversampling parameter $p\geq2$, where $k+p\leq\min\{m,n\}$. Set
$S_k=\sum_{j=k+1}^{\min\{m,n\}}\sigma_j^2$. Draw an $n\times(k + p)$ standard Gaussian matrix $\Omega$ and construct the sample matrix $A\Omega$. Then, for the projection $\Pi_{A\Omega}A$, the expected error satisfies the inequality
$$\mathbb{E}\|A-\Pi_{A\Omega}A\|_F\leq \left(1+\frac{k}{p-1}\right)^{\frac{1}{2}}S_k^{\frac{1}{2}}.$$
\end{Theorem}

We have to be careful when using this result because neither of the two observed approaches completely satisfies the assumptions of Theorem~\ref{tm:error}. Precisely, our Gaussian matrices have $k$ columns, not $k+p$, while the effect of oversampling comes from the fact that the projection matrix $Q$ in Algorithm~\ref{agm:randCMF} contains more than $k$ columns.
However, in order to explain the difference in these two approaches, in this discussion, we will make a minor modification. We are going to assume that $\Omega_1$, $Q_1$, $\Omega_2$, and $Q_2$ from~\eqref{Omega_rand}, as well as $\bar{\Omega}$ and $\bar{Q}$ from~\eqref{Omega12} consist of $k+p$, $p\geq2$, columns.

Then, the approximation of $X$ obtained via $\Omega_1$ and $Q_1$, denoted by $\Pi_{Q_1}X$, satisfies the assumptions of Theorem~\ref{tm:error} and we have
\begin{equation}\label{eq:errorQ1}
\mathbb{E}\|X-\Pi_{Q_1}X\|_F\leq \left(1+\frac{k}{p-1}\right)^{\frac{1}{2}}S_{X,k}^{\frac{1}{2}},
\end{equation}
where $S_{X,k}$ is the sum of the $\min\{m,n_1\}-k$ smallest singular values of $X$. Further on, for the orthogonal matrix $Q$ obtained by the thin QR decomposition, $QR=\begin{bmatrix}
    Q_1 & Q_2
\end{bmatrix}$,
the subspace spanned by the columns of $Q_1$ is a subset of the subspace spanned by the columns of $Q$. Thus, for the projection obtained via $Q$ and denoted by $\Pi_Q$, inequality
\begin{equation}\label{eq:subspace}
\|X-\Pi_QX\|_F\leq\|X-\Pi_{Q_1}X\|_F
\end{equation}
holds. Hence, it follows from the relations~\eqref{eq:errorQ1} and~\eqref{eq:subspace} that
\begin{equation}\label{eq:expX}
\mathbb{E}\|X-\Pi_QX\|_F\leq \left(1+\frac{k}{p-1}\right)^{\frac{1}{2}}S_{X,k}^{\frac{1}{2}},
\end{equation}
In the same way we get
\begin{equation}\label{eq:expY}
\mathbb{E}\|Y-\Pi_QY\|_F\leq \left(1+\frac{k}{p-1}\right)^{\frac{1}{2}}S_{Y,k}^{\frac{1}{2}},
\end{equation}
where $S_{Y,k}$ is the sum of the $\min\{m,n_2\}-k$ smallest singular values of $Y$.
In other words, in the joint approximation obtained via $Q$, approximation error in $X$ depends on its singular values, not on the singular values of $Y$, and vice versa.

On the other hand, if one would take the projection $\Pi_{\bar{Q}}$ obtained via $\bar{Q}$,
Theorem~\eqref{tm:error} implies
\begin{equation}\label{eq:expXY}
\mathbb{E}\|\begin{bmatrix} X & Y \end{bmatrix}-\Pi_{\bar{Q}}\begin{bmatrix} X & Y \end{bmatrix}\|_F\leq \left(1+\frac{k}{p-1}\right)^{\frac{1}{2}}\bar{S}_k^{\frac{1}{2}},
\end{equation}
for $\bar{S}_k=\sum_{j=k+1}^{\min\{m,n_1+n_2\}}\bar{\sigma}_j^2$, where $\bar{\sigma}_j$, $1\leq j\leq\min\{m,n_1+n_2\}$, are the singular values of $\begin{bmatrix} X & Y \end{bmatrix}$. That is, the approximation error in $X$, as well as in $Y$ does not depend directly on the singular values of the matrix in question, but on the singular values of $\begin{bmatrix} X & Y \end{bmatrix}$.

\begin{figure}[h!]
    \begin{minipage}{.5\textwidth}
    \centering
    \includegraphics[width=1\linewidth]{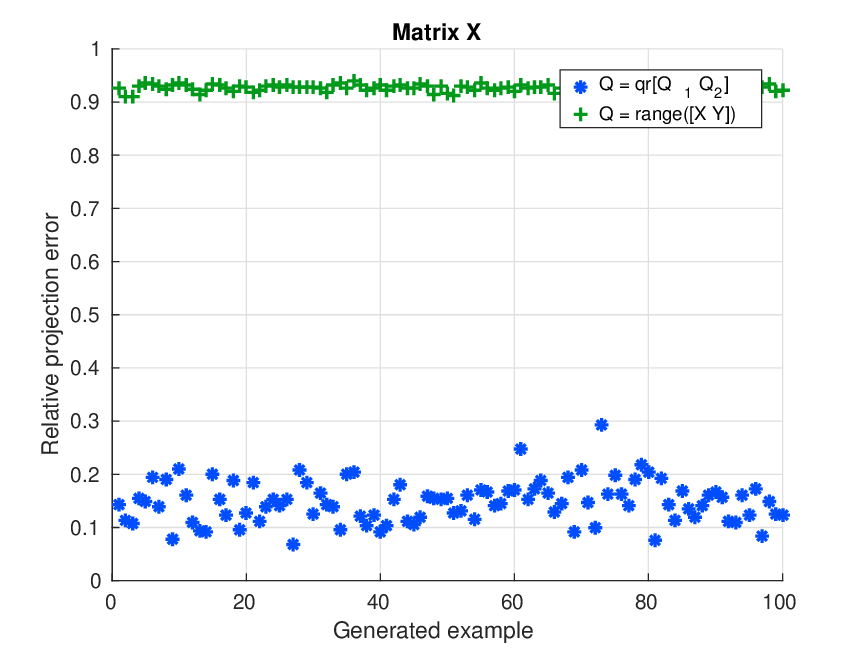}
    \end{minipage}%
    \begin{minipage}{.5\textwidth}
        \centering
        \includegraphics[width=1\textwidth]{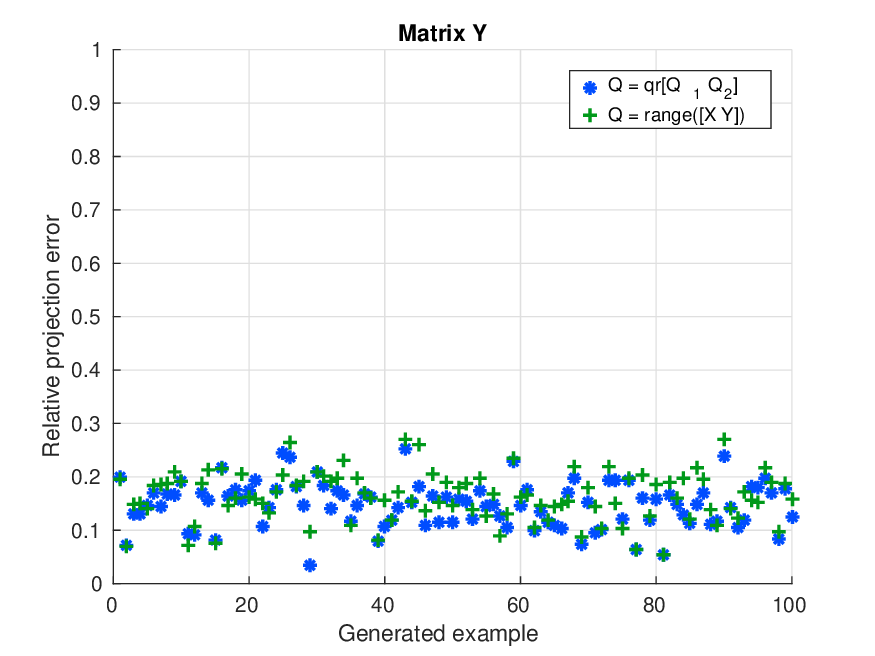}	
    \end{minipage}\\
    \caption{Relative projection errors for matrices $X$ and $Y$ on $100$ randomly generated examples.}
    \label{fig:Projection_err_approx}
\end{figure}

We have observed in our numerical examples that, if the singular values of $X$ and $Y$ are significantly different in size, the approach described by~\eqref{Omega12} results in a much worse relative projection error for the matrix with smaller singular values, compared to the matrix with larger singular values. In contrast, the approach used in Algorithm~\ref{agm:randCMF} results in relative projection errors
that are of the same order of magnitude.
We performed the experiments on $100$ pairs of matrices constructed as follows:
\begin{itemize}
\item \texttt{X = orth(rand($500$,$35$))*diag(sort(rand($35$,$1$)))*(orth(rand($200$,$35$)))';}
\item \texttt{Y = orth(rand($500$,$35$))*diag(100*sort(rand($35$,$1$)))*(orth(rand($300$,$35$)))';}
\end{itemize}
By construction, both matrices have rank $35$. However, the singular values of $\texttt{Y}$ are substantially larger in magnitude. The target rank for coupled decomposition was set to $k=30$. The results given in the relations~\eqref{eq:expX}, \eqref{eq:expY}, and~\eqref{eq:expXY} are illustrated in Figure~\ref{fig:Projection_err_approx}. Our approach achieves a much better projection of the matrix $X$ compared to the method based on $\range(\begin{bmatrix} X & Y \end{bmatrix})$. As demonstrated earlier, this advantage stems from the fact that the singular values of $Y$ are much larger than the singular values of $X$, and they affect the error in $X$. When oversampling is incorporated into the algorithms for randomization, the overall performance improves slightly. However, the relative relationship between the errors remains unchanged. We must emphasize that the results from Figure~\ref{fig:Projection_err_approx} only show that we obtain better projections, they do not guarantee better final approximations.

Additionally, since Algorithm~\ref{agm:randCMF} corresponds to randomized SVD with $\Omega=\begin{bmatrix} \Omega_1 & 0 \\ 0 & \Omega_2\end{bmatrix}\in\R^{(n_1+n_2)\times2k}$, one could find it interesting to compare our randomization approach with randomized SVD of the form~\eqref{Omega12}, but with $\bar{\Omega}$ of size $(n_1+n_2)\times2k$. In that case, the obtained approximation errors are very similar to those obtained by Algorithm~\ref{agm:randCMF}, that is, the gap observed in Figure~\ref{fig:FirstExample} disappears. In order to get matrix $\bar{\Omega}$ one has to perform QR factorization on a matrix with $2k$ columns, while we do two QR factorizations on the matrices with $k$ columns and reorthogonalization of a matrix with $2k$ columns. Although our approach requires more computation in the case of the simple randomization from the Algorithm~\ref{agm:randCMF}, it becomes faster when the advanced randomization techniques are used. As it will be explained in the next subsection, randomized subspace iterations and randomized block Krylov iterations use multiple QR decompositions. In our Algorithms~\ref{agm:rsiCMF} and~\ref{agm:CMFRBKI}, those are the decompositions of the matrices with only $k$, instead of $2k$, columns. Hence, our approach becomes faster, because computational cost of QR decomposition scales quadratically with the number of columns.

\subsection{Refinements of randomized CMF}

In the case of a single matrix, it has been shown that, depending on the singular value distribution of the matrix, more sophisticated randomization techniques can yield greater efficiency. Therefore, we are going to generalize the algorithms for the randomized subspace iteration and the randomized block Krylov iteration from~\cite{TrWe23}, so they can be used with CMF.

The first refinement of Algorithm~\ref{agm:randCMF} uses \emph{randomized subspace iteration (RSI)}. Here, instead of obtaining the orthogonal basis $Q_1$ of $\range(X)$ and $Q_2$ of $\range(Y)$ like in~\eqref{Omega_rand} we have
$$Q_1R_1=(XX^T)^{q-1}X\Omega_1, \quad Q_2R_2=(YY^T)^{q-1}Y\Omega_2,$$
where $\Omega_1\in{\R^{n_1\times k}}$ and $\Omega_2\in{\R^{n_2\times k}}$ are random Gaussian matrices and $q$ is a depth parameter. Typically, $2\leq q\leq5$, see~\cite{TrWe23}.
After we get $Q_1$ and $Q_2$, we proceed in the same way as in Algorithm~\ref{agm:randCMF}.
The procedure for CMF  using RSI is given in Algorithm~\ref{agm:rsiCMF}.

\begin{algorithm}
\caption{RSI CMF}\label{agm:rsiCMF}
\renewcommand{\algorithmicrequire}{\textbf{Input:}}
\renewcommand{\algorithmicensure}{\textbf{Output:}}
\begin{algorithmic}
\Require $X\in\R^{m\times n_1}$, $Y\in\R^{m\times n_2}$, $k<\min\{n_1,n_2\}$
\Ensure $U\in\R^{m\times k}$, $V\in\R^{n_1\times k}$, $W\in\R^{n_2\times k}$
\State Generate random matrices $\Omega_1\in{\R^{n_1\times k}}$ and $\Omega_2\in{\R^{n_2\times k}}$.
\For{$1\leq i\leq q$}
\State $[Q_1,\sim]=\texttt{qr}(X\Omega_1,0)$
\State $\Omega_1=X^TQ_1$
\State $[Q_2,\sim]=\texttt{qr}(Y\Omega_2,0)$
\State $\Omega_2=Y^TQ_2$
\EndFor
\State $[Q,\sim]=\texttt{qr}([Q_1,Q_2],0)$
\State $[\widehat{U},V,W]=\texttt{CMF}(Q^TX,Q^TY,k)$ \Comment{Algorithm~\ref{agm:CMF}}
\State $U=Q\widehat{U}$
\end{algorithmic}
\end{algorithm}

In the \emph{randomized block Krylov iteration} (RBKI) the idea is to project $X$ and $Y$ onto corresponding Krylov subspaces $\mathcal{K}_q(XX^T;X\Omega_1)$ and $\mathcal{K}_q(YY^T;Y\Omega_2)$ defined as
\begin{align*}
    \mathcal{K}_q(XX^T;X\Omega_1) &= \range\begin{bmatrix}
        X\Omega_1 & (XX^T)X\Omega_1 & \ldots & (XX^T)^{q-1}X\Omega_1
    \end{bmatrix},\\
    \mathcal{K}_q(YY^T;X\Omega_2) &= \range\begin{bmatrix}
        Y\Omega_2 & (YY^T)Y\Omega_2 & \ldots & (YY^T)^{q-1}Y\Omega_2
    \end{bmatrix}.
\end{align*}
Matrices $\Omega_1 \in \mathbb{R}^{n_1\times \ell}$ and $\Omega_2\in\mathbb{R}^{n_2\times \ell}$ are random Gaussian matrices. Parameter $\ell$ represents the block dimension, and $q$ is the order of the Krylov subspace. As the result, RBKI method constructs the orthogonal bases $Q_1\in \mathbb{R}^{m\times \ell q}$ of $\mathcal{K}_q(XX^T;X\Omega_1)$ and $Q_2\in \mathbb{R}^{m\times \ell q}$ of $\mathcal{K}_q(YY^T;Y\Omega_2)$. In the literature \cite{HaMaTr11,TrWe23,meyer2024unreasonable}, $\ell$ is usually $1,2,k$ or $k+4$. The parameter $q$ can vary, and we are going to discuss it further within the numerical examples. To define the projections onto $\range (X) + \range(Y)$, we use rank-revealing QR decomposition
and continue the same way as described for Algorithm~\ref{agm:randCMF}. Note that the projection defined for RBKI is usually much larger than the one in the RSI method, due to variability in $\ell$ and $q$. As a consequence, we expect a much more accurate approximation.  We present the algorithm from~\cite{TrWe23} for finding the basis for the Krylov subspaces of order $q$ (Algorithm~\ref{agm:RBKI}), together with our CMF algorithm using RBKI projection (Algorithm~\ref{agm:CMFRBKI}).

\begin{algorithm}
\caption{RBKI \cite{TrWe23}}\label{agm:RBKI}
\renewcommand{\algorithmicrequire}{\textbf{Input:}}
\renewcommand{\algorithmicensure}{\textbf{Output:}}
\begin{algorithmic}
\Require $A\in\R^{m\times n}$, block size $\ell$, iteration count $q$
\Ensure Matrix $Q\in \R^{m\times \ell q}$ representing basis for Krylov subspace $\mathcal{K}_q(AA^T;A\Omega_0)$
\State Generate random $\Omega_0 \in \R^{n \times \ell}$
\For{$i=1,\ldots,q$}
\State $Q_i = A\Omega_{i-1}$
\State $Q_i = Q_i - \sum_{j<i}Q_j(Q_j^TQ_i)$
\State $Q_i = Q_i - \sum_{j<i}Q_j(Q_j^TQ_i)$ \Comment{Reorthogonalization}
\State $[Q_i,\sim] = $ \texttt{qr}$(Q_i,0)$
\State $\Omega_i = A^TQ_i$
\EndFor
\State $Q = [Q_1 \ldots Q_q]$
\end{algorithmic}
\end{algorithm}

\begin{algorithm}
\caption{RBKI CMF}\label{agm:CMFRBKI}
\renewcommand{\algorithmicrequire}{\textbf{Input:}}
\renewcommand{\algorithmicensure}{\textbf{Output:}}
\begin{algorithmic}
\Require $X\in\R^{m\times n_1}$, $Y\in\R^{m\times n_2}$, $\ell$, $q$, $k < \min\{n_1,n_2\}$
\Ensure $U\in \R^{m\times k}$, $V\in \R^{n_1\times k}$, $W\in \R^{n_2\times k}$
\State $Q_1 = \texttt{RBKI}(X,\ell,q)$ \Comment{Algorithm~\ref{agm:RBKI}}
\State $Q_2 = \texttt{RBKI}(Y,\ell,q)$
\State $[Q,\sim] = $\texttt{qr}$([Q_1, Q_2],0)$ \Comment{$Q\in\R^{n\times 2kq}$}
\State $[U,V,W] = $\texttt{CMF}$(Q^TX, Q^TY, k)$ \Comment{Algorithm~\ref{agm:CMF}}
\State $U = QU$
\end{algorithmic}
\end{algorithm}

The comparison of Algorithms~\ref{agm:CMF}, \ref{agm:randCMF}, \ref{agm:rsiCMF}, and~\ref{agm:CMFRBKI} through the numerical results will be presented in the next section.

\section{Numerical examples for CMF}\label{sec:CMF_num}

In this section, we construct numerical experiments to test the efficiency of our algorithms and to compare proposed randomized methods, i.e., Algorithm~\ref{agm:randCMF}, Algorithm~\ref{agm:rsiCMF}, and Algorithm~\ref{agm:CMFRBKI}. We use MATLAB 9.0.0.341360 (R2016a) in double precision (IEEE Standard 754).

The accuracy of our algorithms is determined by comparing the corresponding relative errors,
\begin{equation}\label{eq:err}
err_X = \frac{\|X - UV^T\|_F}{\|X\|_F}, \quad err_Y = \frac{\|Y - UW^T\|_F}{\|Y\|_F}.
\end{equation}
The results are presented in the following way. The relative error for the basic algorithm (Algorithm~\ref{agm:CMF}) serves as the benchmark. We vary parameters for both the RSI and RBKI algorithms.
In our tables, $k$ is the approximation rank, $p$ is the oversampling parameter, $\ell$ is the block dimension for RBKI, and $q$ is the depth parameter in RSI, that is, the order of the Krylov subspace in RBKI.
For RSI, we report one selected result.
We test the RBKI algorithm using the blocks of dimension $\ell=1,2,k$ and $k+4$. This results in different oversampling parameters $p$.
For each block dimension, we highlight one result that is close enough to the benchmark result, along with the corresponding oversampling parameter. Furthermore, for selected problems, for the block dimensions $\ell=1$ and $\ell=2$, we present the plots showing the relative errors as a function of the oversampling parameter. These plots demonstrate that even with a smaller oversampling parameter than the one shown in the table, the error remains close to that of the basic algorithm.

Let us present the test matrices and then analyze each example individually.
\begin{enumerate}
    \item[1.] \texttt{SyntheticTest1}: We construct low rank matrices $X\in\R^{m\times n_1}, Y\in\R^{m\times n_2}$ with no special structure for singular values as in Section \ref{ssec:randCMF}:
    \begin{itemize}
        \item \texttt{X = rand(m,$r_1$)*rand($r_1$,$n_1$);}
        \item \texttt{Y = rand(m,$r_2$)*rand($r_2$,$n_2$);}
    \end{itemize}

    \item [2.] \texttt{SyntheticTest2}: Matrix $X\in \R^{n\times n}$ is constructed as in~\cite{saibaba2023randomized}, i.e., for $d\in\mathbb{N}$:
    \begin{itemize}
        \item $\Sigma_X = \operatorname{diag}( {1,\ldots,1},2^{-d},3^{-d},\ldots,(n-r+1)^{-d})$;
        \item \texttt{UX = orth(rand(n));} \texttt{VX = orth(rand(n));}
        \item \texttt{X = UX*$\Sigma_X$*VX';}
    \end{itemize}
    Notice that $X$ is designed to have quickly decaying singular values. Next, the matrix $Y\in \mathbb{R}^{n\times n}$ is constructed so that parts of $\range(X)$ and $\range(Y)$ intersect. This is controlled with the parameter $c$:
    \begin{itemize}
        \item $\Sigma_Y = \operatorname{diag}(1,\ldots,1,2^{-1},3^{-1},\ldots,(n-2r+1)^{-1})$;
        \item \texttt{UY = ([UX(:,1:c) orth(rand(n,n-c))];}
        \item \texttt{VY = ([VX(:,1:c) orth(rand(n,n-c))];}
        \item \texttt{Y = UY*$\Sigma_Y$*VY';}
    \end{itemize}

    \item[3.] \texttt{SynthethicTest3}: Matrix $X\in \R^{m\times n}$ is constructed as in \cite{saibaba2023randomized}. More precisely,
    \begin{equation*}
        X = \sum^r_{j=1}\frac{10}{j}x^X_j\left(y^X_j\right)^T + \sum^{\min\{m,n\}}_{j=r+1} \frac{1}{j}x^X_j\left(y^X_j\right)^T,
    \end{equation*}
    where $x^X_j\in\R^{m}$ and $y^X_j\in\R^n$ are sparse random vectors with density $0.25$, created using MATLAB command \texttt{sprand}.
    We construct the matrix $Y\in\R^{m\times n}$ in a similar way, with the first $r$ random vectors $x^X_j\in\R^m$, $y^X_j\in\R^n$ the same as for $X$,
    \begin{equation*}
        Y = \sum^r_{j=1}\frac{10}{j}x^X_j\left(y^X_j\right)^T + \sum^{\min\{m,n\}}_{j=r+1} \frac{1}{j}x^Y_j\left(y^Y_j\right)^T.
    \end{equation*}

    \item [4.] \texttt{SyntheticTest4}: Matrix $X\in\R^{m\times n_1}$ is defined as:
    \begin{itemize}
        \item $\Sigma_X = \operatorname{diag}(2^{-i})_{i=1}^{n_1};$
        \item \texttt{X=orth(rand(m))*$\Sigma_X$;}
    \end{itemize}
    while $Y\in\R^{m\times n_2}$ is a low-rank matrix with no special structure:
    \begin{itemize}
        \item \texttt{Y = rand(m,$r_2$)*rand($r_2$,$n_2$);}
    \end{itemize}

    \item [5.]\texttt{SyntheticTest5}: Here, both $X\in\R^{m\times n_1}$ and $Y\in\R^{m\times n_2}$ are ill-conditioned with fast decaying singular values:
    \begin{itemize}
        \item \texttt{$\mathtt{A} = \operatorname{diag}(1, 2^{-1}, 2^{-2},\ldots, 2^{-n_1})$};
        \item \texttt{UA = orth(rand(m))};
        \item \texttt{A = UA*A};
        \item \texttt{$\mathtt{B} = \operatorname{diag}(1, 2^{-1}, 2^{-2},\ldots, 2^{-n_2})$};
        \item \texttt{UB = [UA(1:10) orth(rand(m,m-10))]};
        \item \texttt{B = UB*B};
    \end{itemize}
\end{enumerate}

\subsection{SyntheticTest1}
\begin{figure}[ht]
    \centering
    \begin{minipage}{.5\textwidth}
\includegraphics[width=1\textwidth]{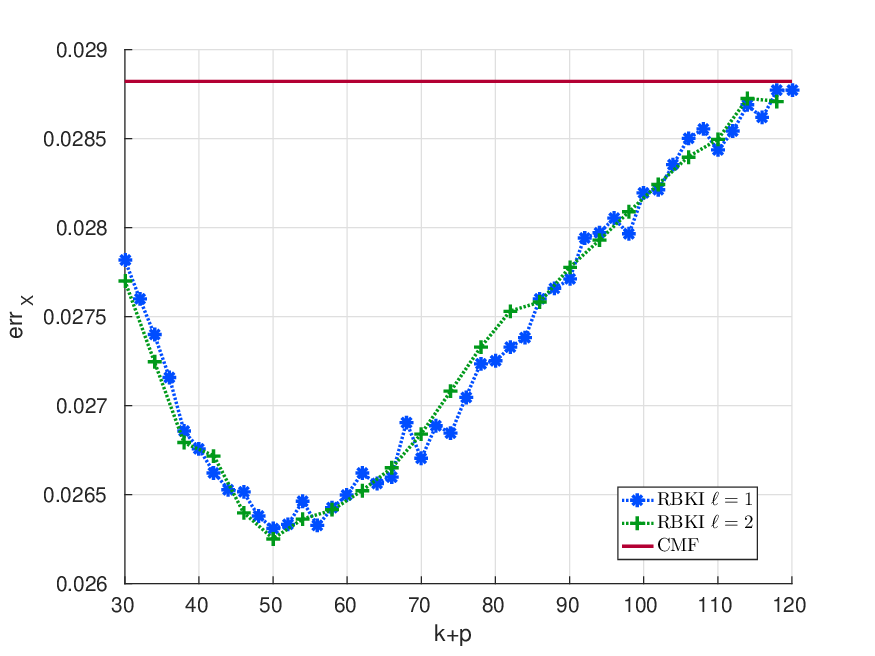}
    \end{minipage}%
    \begin{minipage}{.5\textwidth}
        \centering
\includegraphics[width=1\textwidth]{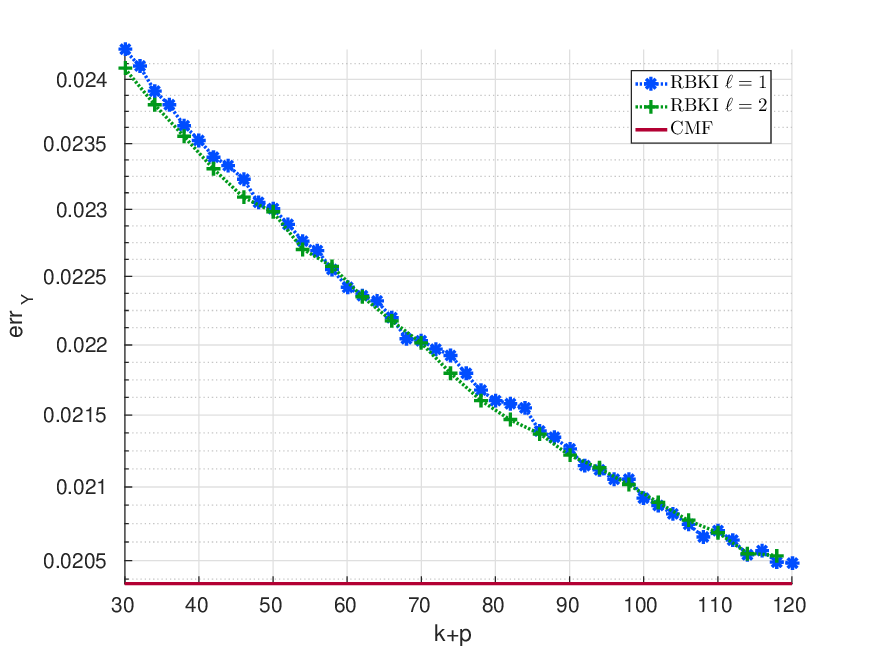}	
    \end{minipage}
    \caption{Relative errors with the respect to the dimension of the projection subspace in Algorithm~\ref{agm:CMFRBKI} for \texttt{SyntheticTest1}.}
    \label{fig:diff4}
\end{figure}

We consider \texttt{SyntheticTest1} with $m=500$, $n_1=200$, $n_2=300$, $r_1=100$ and $r_2=150$. We are looking for the low-rank approximation of order $k=30$.
%ovdje samo povezati tekst
We present the performance of Algorithms~\ref{agm:rsiCMF} and~\ref{agm:CMFRBKI}.
The selected results are given in Table~\ref{tab:testdiff4}. Figure~\ref{fig:diff4} shows the relative errors for different oversampling parameters in RBKI. Notice that the relative error for $X$ initially decreases below the benchmark, then increases until it reaches the benchmark error. In contrast, the error for $Y$ consistently decreases. We assume that this happens because $X$ and $Y$ do not share any common features and it takes a number of iterations for the relative errors to balance.

\begin{table}[ht]
        \centering
        \begin{tabular}{|l|c|c|c|c|c|}
            \hline
            Algorithm & $p$ & $\ell$ & $q$ & $err_X$ & $err_Y$ \\ \hline\hline
            Basic CMF & - & - & - & $2.88218893\cdot 10^{-2}$ & $2.03454946 \cdot 10^{-2}$ \\
            RSI & 30 & - & 5 & $2.86340153\cdot 10^{-2}$ & $2.06624393\cdot 10^{-2} $ \\
            RBKI & 112 & 1 & 71 & $2.88153938\cdot 10^{-2}$ & $2.04048432\cdot 10^{-2}$\\
            RBKI & 122 & 2 & 38 & $2.89287460\cdot 10^{-2}$ & $ 2.03281104 \cdot 10^{-2}$\\
            RBKI & 90 & 30 & 2 & $2.83115828\cdot 10^{-2}$ & $2.11079710\cdot 10^{-2}$\\
            RBKI & 106 & 34 & 3 & $2.83986717\cdot 10^{-2}$ & $2.09275944\cdot 10^{-2}$\\
             \hline
        \end{tabular}
        \smallskip
        \caption{Selected result for \texttt{SyntheticTest1} and $k=30$. Oversampling parameters and relative errors for RSI with $q=5$ and RBKI with $\ell=1,2,k,k+4$.}
        \label{tab:testdiff4}
\end{table}

\subsection{SyntheticTest2}
We consider $n=1000$, $r=15$, $d=2$ and $c=50$. We are looking for the low-rank approximation of order $k=50$.

Selected results for Algorithms~\ref{agm:randCMF}, \ref{agm:rsiCMF}, and~\ref{agm:CMFRBKI} are presented in Table~\ref{tab:testdiff2}. We also give the errors obtained by the basic algorithm for the comparison.
For this and the following example, we document the running time for each algorithm using MATLAB function \texttt{tic toc}. For Algorithms~\ref{agm:randCMF}, \ref{agm:rsiCMF}, and \ref{agm:CMFRBKI} we present the total time, i.e., the time needed for constructing the projection matrix and running CMF on a smaller set of matrices. Separately, we present the time for the CMF algorithm itself. For the basic CMF total time and CMF time are the same, since there is no preprocessing, the matrices are not projected. We can see that, as expected, the randomized algorithms are always faster than the basic CMF.

\begin{figure}[ht]
    \centering
    \begin{minipage}{.5\textwidth}
        \includegraphics[width=1\textwidth]{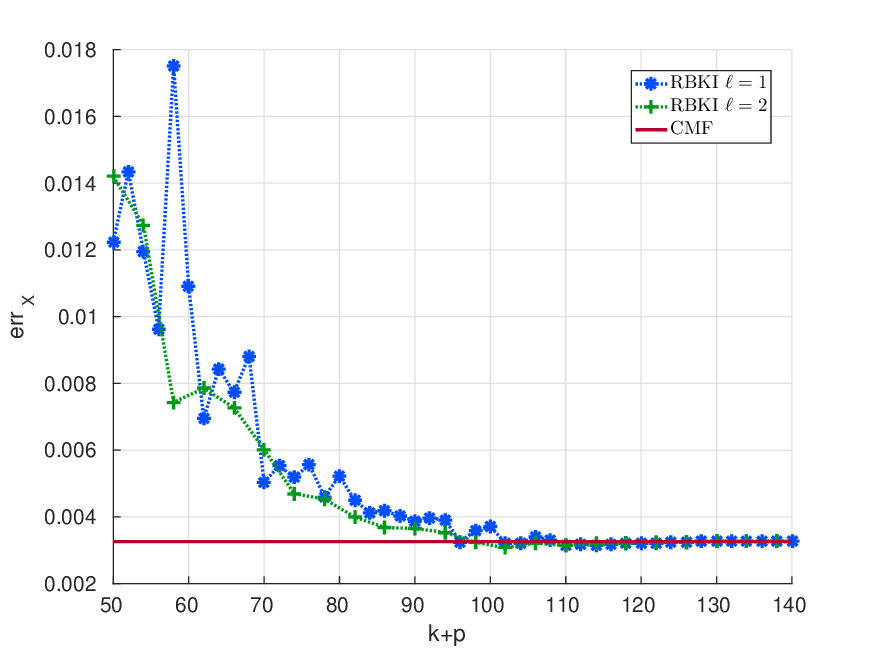}
    \end{minipage}%
    \begin{minipage}{.5\textwidth}
        \centering
	\includegraphics[width=1\textwidth]{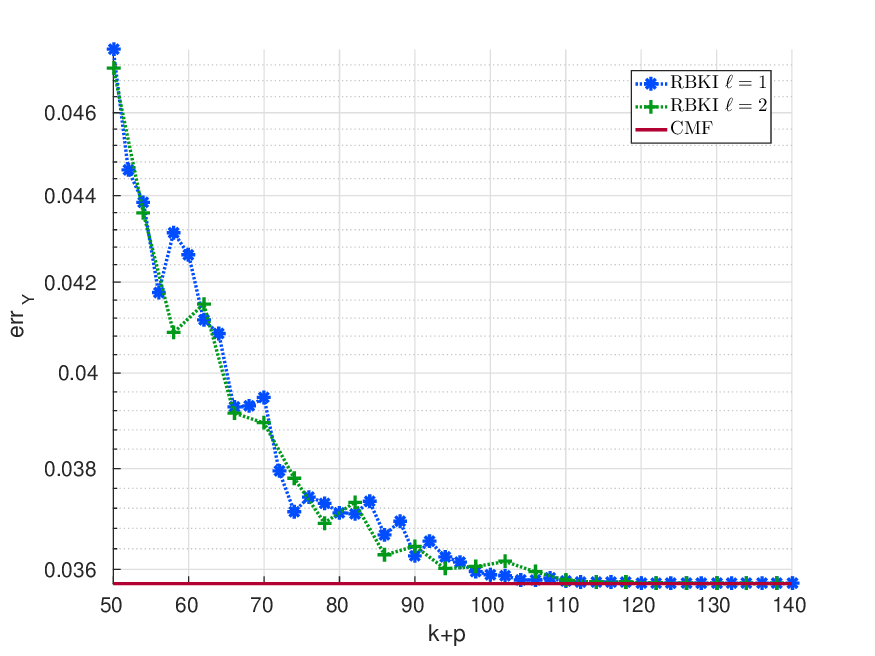}	
    \end{minipage}
    \caption{Relative errors with the respect to the dimension of the projection subspace in Algorithm~\ref{agm:CMFRBKI} for \texttt{SyntheticTest2}.}
    \label{fig:diff2}
\end{figure}

The oversampling parameter for RSI is $p=k=50$, but we need $q=4$ iterations to bring the relative error close to the one of the basic algorithm. The errors are slightly better for the RBKI algorithm. The oversampling parameters for RBKI are higher because of the way we construct the projection matrices. However, even a smaller oversampling parameter provides a good enough approximation, which is presented in Figure~\ref{fig:diff2}. Simple randomization, i.e., Algorithm~\ref{agm:randCMF}, performed the worst, which is anticipated.

\begin{table}[ht]
        \centering
        \begin{tabular}{|l|c|c|c|c|c|c|c|}
            \hline
            Algorithm  & $p$ & $\ell$ & $q$ & $err_X$ & $err_Y$ & \makecell{total\\ time (s)} & \makecell{CMF\\ time (s)} \\ \hline\hline
            Basic CMF  & - & - & - & ${3.25997151\cdot 10^{-3}}$ & ${3.57247002 \cdot 10^{-2}}$ & ${0.288542}$ & ${0.288542}$ \\
            Randomized  & 50 & - & - & $2.13494569 \cdot 10^{-3}$ & $4.19216919\cdot 10^{-2}$ & $0.017870$ & $0.009305$\\
            RSI  & 50 & - & 4 & $3.25992287\cdot 10^{-3}$ & $ 3.57585534\cdot 10^{-2} $ & $0.025313$ & $0.008786$ \\
            RBKI  & 84 & 1 & 67 & $3.25997444\cdot 10^{-3}$ & $3.57247002\cdot 10^{-2}$ & $0.125757$ & $0.012734$\\
            RBKI  & 86 & 2 & 34 & $3.25997518\cdot 10^{-3}$ & $ 3.57247006 \cdot 10^{-2}$ & $0.077220$ & $0.012798$\\
            RBKI  & 150 & 50 & 2 & $3.25979677\cdot 10^{-3}$ & $3.58010584\cdot 10^{-2}$ & $0.044432$ & $0.020858$\\
            RBKI  & 166 & 54 & 2 & $3.25980127\cdot 10^{-3}$ & $3.57555301\cdot 10^{-2}$ & $0.050801$ & $0.024738$\\
            \hline
        \end{tabular}
        \smallskip
        \caption{Selected results for \texttt{SyntheticTest2} and $k=50$. Oversampling parameters, relative errors for RSI with $q=4$ and RBKI with $\ell=1,2,k,k+4$, and running times.}
        \label{tab:testdiff2}
    \end{table}

\subsection{SyntheticTest3}
For this experiment, we consider $m=10000$, $n=500$ and $r=50$. We are looking for the rank $k=30$ approximation.

Selected results are presented in Table~\ref{tab:testdiff3}.
Regarding the execution time, the difference between the randomized algorithms and Algorithm~\ref{agm:CMF} is noticeable. This is because the randomized algorithms carry out CMF of the much smaller matrices.
For the RSI algorithm, we needed $q=5$ iterations to match the benchmark error. The oversampling parameter in the RBKI method is greater. However, by examining the result, we can conclude that, like in the previous example, even a smaller oversampling parameter leads to a satisfactory result. We can also conclude that a small $\ell$ in RBKI gives as good result as a larger $\ell$. An advantage of a smaller block size is a smaller oversampling parameter, resulting in the smaller matrices on which CMF is performed. However, the total execution time is shorter for the larger blocks because the RBKI algorithm demands only two iterations. Again, Algorithm~\ref{agm:randCMF} has the worst performance.

\begin{table}[ht]
        \centering
        \begin{tabular}{|l|c|c|c|c|c|c|c|}
            \hline
            Algorithm  & $p$ & $\ell$ & $q$ & $err_X$ & $err_Y$ & \makecell{total\\ time (s)} & \makecell{CMF\\ time (s)} \\ \hline\hline
            Basic CMF  & - & - & - & $6.03560458\cdot 10^{-2}$ & $6.03576711 \cdot 10^{-2}$ & $0.527464$ & $0.527464$\\
            Randomized  & 30 & - & - & $6.15615406\cdot 10^{-2}$ & $6.15000752\cdot 10^{-2}$ & $0.054951$ & $0.008305$\\
            RSI  & 30 & - & 5 & $6.03862191\cdot 10^{-2}$ & $ 6.03867712\cdot 10^{-2} $ & $0.127206$ & $0.008582$ \\
            RBKI  & 64 & 1 & 47 & $6.03560471\cdot 10^{-2}$ & $6.03576722\cdot 10^{-2}$ & $0.562032$ & $0.006410$\\
            RBKI  & 74 & 2 & 26 & $6.03560459\cdot 10^{-2}$ & $6.03576711 \cdot 10^{-2}$ & $0.426180$ & $0.006848$\\
            RBKI  & 90 & 30 & 2 & $6.03685901\cdot 10^{-2}$ & $6.03696584\cdot 10^{-2}$ & $0.146570$ & $0.008841$\\
            RBKI  & 106 & 34 & 2 & $6.03610827\cdot 10^{-2}$ & $6.03625993\cdot 10^{-2}$ & $0.127120$ & $0.009377$\\
             \hline
        \end{tabular}
                \smallskip
        \caption{Selected results for \texttt{SyntheticTest3}  and $k=30$. Oversampling parameters, relative errors for RSI with $q=5$ and RBKI with $\ell=1,2,k,k+4$, and running times.}
        \label{tab:testdiff3}
\end{table}

\subsection{SyntheticTest4}
Consider $m=500$, $n_1 =300$, $n_2=200$ and $r_2=100$. We want to find the rank $k=30$ approximation.

Table \ref{tab:testdiff6} represents the selected results. Similarly as in \texttt{SynthethicTest1}, matrices $X$ and $Y$ do not have anything in common, thus, the behavior is very much alike. The difference here is that it takes a larger oversampling parameter to obtain a good error for $Y$.

\begin{table}[ht]
        \centering
        \begin{tabular}{|l|c|c|c|c|c|}
            \hline
            Algorithm  & $p$ & $\ell$ & $q$ & $err_X$ & $err_Y$ \\ \hline\hline
            Basic CMF  & - & - & - & $4.86320647\cdot 10^{-1}$ & $2.02998694 \cdot 10^{-2}$\\
            RSI  & 30 & - & 8 & $4.86654497\cdot 10^{-1}$ & $ 2.04024796\cdot 10^{-2} $ \\
            RBKI  & 118 & 1 & 74 & $4.86320519\cdot 10^{-1}$ & $2.02998694\cdot 10^{-2}$\\
            RBKI  & 118 & 2 & 37 & $4.86315021\cdot 10^{-1}$ & $ 2.02998748 \cdot 10^{-2}$\\
            RBKI  & 60 & 30 & 2 & $4.85925751\cdot 10^{-1}$ & $2.09976041\cdot 10^{-2}$\\
            RBKI  & 68 & 34 & 2 & $4.85065718\cdot 10^{-1}$ & $2.08038065\cdot 10^{-2}$\\
             \hline
        \end{tabular}
        \smallskip
        \caption{Selected results for \texttt{SyntheticTest4} and $k=30$. Oversampling parameters and relative errors for RSI with $q=8$ and RBKI with $\ell=1,2,k,k+4$, and running times.}
        \label{tab:testdiff6}
    \end{table}

\subsection{SyntheticTest5. Range intersection.}
This example is constructed to demonstrate the importance of using rank-revealing QR decomposition when constructing the projection matrix $Q$ in the randomized algorithms. This results in a possibly smaller oversampling parameter $p$, as discussed in the Subsection~\ref{ssec:randCMF}.

We consider $m = 500, n_1 = 300$, $n_2 = 200$, and we want to find the approximation of rank $k=30$. We separately analyze the influence of rank determination when using RBKI with block dimensions $\ell=1$ and $\ell=2$.

\begin{figure}[ht]
    \centering
    \begin{minipage}{.45\textwidth}
        \centering
	\includegraphics[width=1\textwidth]{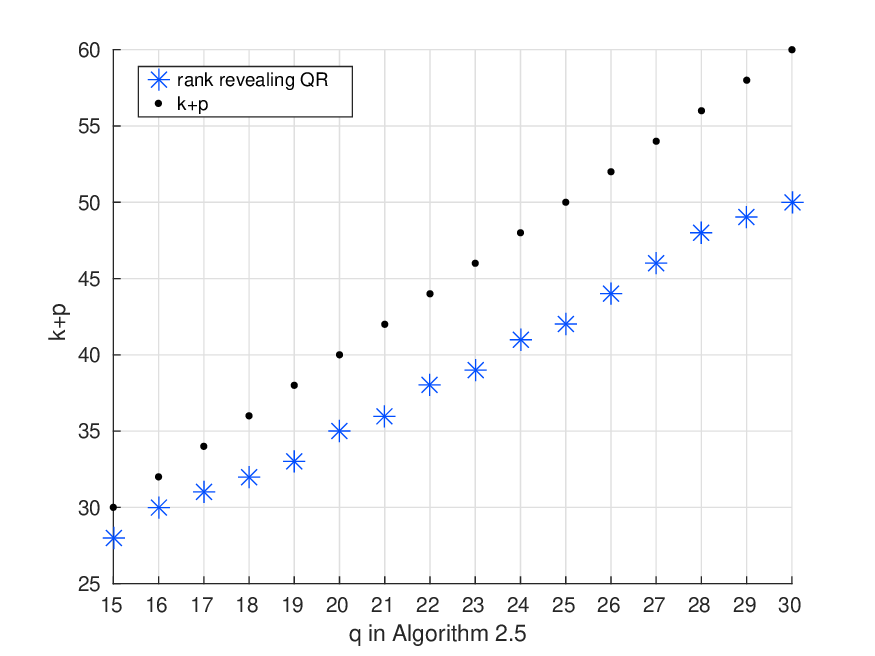}	
    \end{minipage}
    \begin{minipage}{.45\textwidth}
        \centering
	\includegraphics[width=1\textwidth]{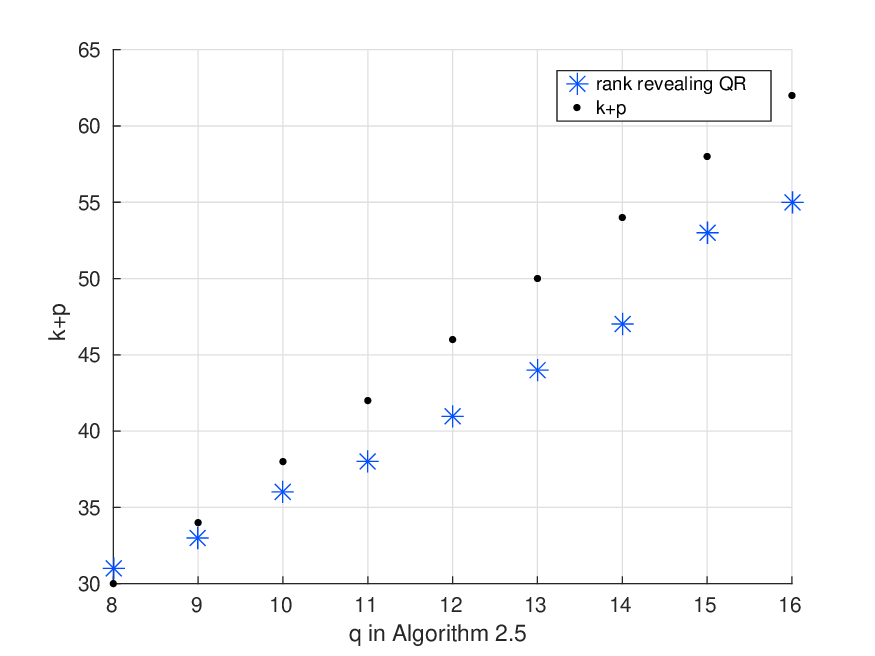}	
    \end{minipage}
    \caption{RBKI with $\ell=1$ (on the left) and $\ell=2$ (on the right). Rank-revealing QR leads to a smaller oversampling parameter.}
    \label{fig:RBKI5k2}
\end{figure}

For {$\ell = 1$} we iterate $q$ from $15$ to $30$ to get the oversampling parameter $p\leq (2q\ell-k)$. For $\ell=2$, we iterate $q$ from 8 to 16 to get the oversampling parameter $p\leq (4q\ell-k)$. Figure \ref{fig:RBKI5k2} demonstrates the difference between the rank of the matrix $Q$ and the maximal dimension of the matrix $Q$.
The impact of the rank-revealing QR should be more noticeable on large-dimensional matrices. Randomization itself contributes to the reduction of dimensionality and the CMF algorithm is executed on smaller matrices. QR with pivoting can further reduce that dimensionality.

Let us end this section with a remark on another well-known randomization technique --- Generalized Nystrom (GN), \cite{nakatsukasa2020fast}.

\begin{Remark}
Generalized Nystrom is a method that produces a rank-$k$ approximation of any matrix $X \in\mathbb{R}^{m\times n}$ using random matrices. We combined stabilized GN with CMF and tested it on a few examples. The results were compared with other methods presented in this paper. The comparison showed that CMF with GN is slower and less accurate. Other methods are faster because they first approximate matrices $X$ and $Y$ with $QQ^{T}X$ and $QQ^{T}Y$, respectively. By determining the matrix $Q$, we have already found one part of the matrix $U$ (which is a common part of the matrices X and Y). What remains is to find the CMF of the smaller matrices $Q^{T}X$ and $Q^{T}Y$. On the other hand, GN does not have this feature, which makes its combination with CMF slower. In addition, all methods were tested with Gaussian random matrices, which in the case of GN produces high errors. The results improved for GN when we replaced those matrices with a random subset of the columns of the identity matrix, but they were still less accurate when compared to the other methods. There is a possibility that more carefully chosen random matrices would produce more accurate results, but the question of algorithm speed remains.
\end{Remark}

\section{Coupled matrix and tensor factorization}\label{sec:CMTF}

Coupled matrix and tensor factorization (CMTF) can be approached differently based on the tensor format. Here we consider two: the Tucker format in Subsection~\ref{ssec:CMTF_Tucker} and the CP format in Subsection~\ref{ssec:CMTF_CP}. For the sake of simplicity, we assume that the tensors are of order three. Generalization on order-$d$, $d>3$, tensors is direct but requires more complicated computations.

\subsection{Tensor preliminaries}

Throughout the text, we denote the tensors by calligraphic letters, e.g., $\calX$. The \emph{order} of a tensor is its dimension. Thus, $\calX\in\R^{n_1\times n_2\times n_3}$ is an order-$3$ tensor.
Vectors obtained from a tensor by fixing all indices but the $n \textsuperscript{th}$ one are called \emph{mode-$n$ fibers}. \emph{Mode-$n$ matricization} (unfolding) is a matrix representation of a tensor acquired by arranging mode-$n$ fibers into a matrix. Mode-$n$ matricization of a tensor $\calX$ is denoted by $X_{(n)}$. Mode-$1$ matricization of $\calX\in\R^{n_1\times n_2\times n_3}$ is an $n_1\times n_2n_3$ matrix. \emph{Mode-$n$ product} of a tensor $\calX$ and a matrix $M$ is a tensor $\calT$,
$$\calT=\calX\times_n M,$$
such that
\begin{equation}\label{product}
T_{(n)}=MX_{(n)}.
\end{equation}
For any matrices $M_1$ and $M_2$ of the appropriate size, equality
\begin{equation}\label{productmn}
\calX\times_mM_1\times_nM_2=\calX\times_nM_2\times_mM_1, \quad \text{if } m\neq n,
\end{equation}
holds.
\emph{Tensor norm} is a generalization of the Frobenius matrix norm,
$$\|\calX\|=\sqrt{\sum_{i=1}^{n_1}\sum_{j=1}^{n_2}\sum_{k=1}^{n_3} x_{ijk}^2}, \quad \calX\in\R^{n_1\times n_2\times n_3}.$$

For this paper, we need two tensor decompositions. The first one is the \emph{Tucker decomposition}~\cite{KoBa09,Tucker66}. It is a decomposition of a tensor into a core tensor multiplied by a matrix in each mode. The Tucker decomposition of $\calX$ takes the form
$$\calX=\calS\times_1M_1\times_2M_2\times_3M_3.$$
\emph{The multilinear rank} of $\calX$ is a triplet $(r_1,r_2,r_3)$, where $r_i=\rank(X_{(i)})$, $i=1,2,3$. If
\begin{equation}\label{multilinrank}
\calX\approx\calS\times_1V_1\times_2V_2\times_3V_3
\end{equation}
and $\calS\in\R^{r_1\times r_2\times r_3}$, we say that~\eqref{multilinrank} is a multilinear rank-$(r_1,r_2,r_3)$ approximation of $\calX$.

The other decomposition of our interest is \emph{CP decomposition}~\cite{KoBa09,KiersCP}. It is a decomposition of a tensor $\calX$ into a sum of rank-one tensors. A rank-one order-$3$ tensor is a tensor that can be written as an outer product of three vectors,
$$\calT=a\circ b\circ c,$$
where $\circ$ stands for the outer product. Then, we write
\begin{equation}\label{exactCP}
\calX=\sum_{i=1}^R a_i\circ b_i\circ c_i\equiv[[A,B,C]],
\end{equation}
for $A=\begin{bmatrix}
a_1 & a_2 & \cdots & a_R
\end{bmatrix}$, $B=\begin{bmatrix}
b_1 & b_2 & \cdots & b_R
\end{bmatrix}$, $C=\begin{bmatrix}
c_1 & c_2 & \cdots & c_R
\end{bmatrix}$.
\emph{Tensor rank} is the smallest number $R$ in the CP decomposition~\eqref{exactCP} and
\begin{equation}\label{trank}
\calX\approx\sum_{i=1}^r \tilde{a}_i\circ \tilde{b}_i\circ \tilde{c}_i\equiv[[\tilde{A},\tilde{B},\tilde{C}]]
\end{equation}
is a rank-$r$ approximation of $\calX$.

Before we proceed with CMTF, let us also define two matrix products that will be used in this section. The \emph{Khatri-Rao product}, denoted by $\odot$, is a product of two matrices $A\in\R^{m\times n}, B\in \R^{p\times n}$,
$$A\odot B = \begin{bmatrix}
a_1 \otimes b_1 & a_2 \otimes b_2 & \cdots & a_n \otimes b_n
\end{bmatrix}
\in\R^{(mp)\times n},$$
where $a_k$ and $b_k$ denote the $k$th column of $A$ and $B$, respectively.
An important property of CP decomposition is that, if~\eqref{exactCP} holds, then
\begin{equation}\label{CPmat}
X_{(1)}=A(C\odot B)^T, \quad X_{(2)}=B(C\odot A)^T, \quad X_{(3)}=C(B\odot A)^T.
\end{equation}
Moreover, the \emph{Hadamard product}, denoted by $\ast$, of two matrices $A,B\in\R^{m\times n}$, is their element-wise product,
$$A\ast B=\begin{bmatrix}
a_{11}b_{11} & a_{12}b_{12} & \cdots &a_{1n}b_{1n} \\
a_{21}b_{21} & a_{22}b_{22} & \cdots &a_{2n}b_{2n} \\
\vdots & \vdots & \ddots & \vdots\\
a_{m1}b_{m1} & a_{m2}b_{m2} & \cdots &a_{mn}b_{mn}
\end{bmatrix}
\in\R^{m\times n}.$$

\subsection{Tensor in Tucker format}\label{ssec:CMTF_Tucker}

Let $\calX\in\R^{m\times n_2\times n_3}$ and $Y\in\R^{m\times n}$. Joint rank-$k$ approximation of $\calX$ and $Y$, coupled in the first mode and taking $\calX$ in the Tucker format, as it is shown in Figure~\ref{fig:CMTFTucker}, is given by
\begin{equation}\label{CMTF_Tucker}
\calX\approx \calS\times_1U\times_2V_2\times_3V_3, \quad Y\approx UW^T,
\end{equation}
where $U\in\R^{m\times k}$, $V_2\in\R^{n_2\times k}$, $V_3\in\R^{n_3\times k}$, $W\in\R^{n\times k}$, $\calS\in\R^{k\times k\times k}$, and multilinear rank of $\calX$ equals $(k,k,k)$.

\begin{figure}[t]
\begin{tikzpicture}
\draw [draw=black] (0,0) rectangle (2,2);
\node at (1,1) {$\calX$};
\draw (0,2)--(1,2.5);
\draw (2,2)--(3,2.5);
\draw (2,0)--(3,0.5);
\draw (1,2.5)--(3,2.5);
\draw (3,0.5)--(3,2.5);
\node at (3.5,1) {$\approx$};
\end{tikzpicture}
\begin{tikzpicture}
\draw [fill=gray, draw=gray] (0,0) rectangle (0.2,2);
\draw [draw=black] (0.4,1.8) rectangle (0.6,2);
\draw (0.4,2)--(0.5,2.1);
\draw (0.6,2)--(0.7,2.1);
\draw (0.6,1.8)--(0.7,1.9);
\draw (0.5,2.1)--(0.7,2.1);
\draw (0.7,2.1)--(0.7,1.9);
\draw [draw=black] (0.9,1.9) rectangle (2.9,2.1);
\draw (0.5,2.3)--(0.7,2.3);
\draw (0.9,2.8)--(1.1,2.8);
\draw (0.5,2.3)--(0.9,2.8);
\draw (0.7,2.3)--(1.1,2.8);
\node at (0.5,0) {$U$};
\end{tikzpicture}
\caption{Graphical description of the tensor $\calX$ from CMTF~\eqref{CMTF_Tucker}.}
\label{fig:CMTFTucker}
\end{figure}

Using the properties of the mode-$n$ product~\eqref{product} and~\eqref{productmn}, tensor $\calX$ from~\eqref{CMTF_Tucker} can be represented as
\begin{equation}\label{CMTF_TuckerU}
X_{(1)}\approx U(\calS\times_2V_2\times_3V_3)_{(1)}.
\end{equation}
Thus, to get a solution of CMTF~\eqref{CMTF_Tucker} we can observe CMF of the tensor matricization $X_{(1)}$ and matrix $Y$. We get the approximations
\begin{equation}\label{CMTF_Tucker2}
X_{(1)}\approx UV^T, \quad Y\approx UW^T.
\end{equation}
The problem of finding the matrices $U$, $V$, and $W$ from~\eqref{CMTF_Tucker2} is the same as the problem~\eqref{CMF} and Algorithm~\ref{agm:CMF} can be used.
Now, the approximation $\widetilde{X}$ can be formed by folding $UV^T$ into a tensor. The procedure is given in Algorithm~\ref{agm:CMTF-Tucker}.

If there are $\calS$, $V_2$ and $V_3$ of the appropriate size, such that $V\in\R^{n_2n_3\times k}$ can be written as $V^T=(\calS\times_2V_2\times_3V_3)_{(1)}$, approximation reached in Algorithm~\ref{agm:CMTF-Tucker} would be the best possible rank-$k$ approximation of the form~\eqref{CMTF_Tucker}.
Although this cannot be guaranteed, a simple procedure in Algorithm~\ref{agm:CMTF-Tucker} extracts the joint factor $U$, and returns a competitive approximation. Its application is presented in Section~\ref{sec:facerec}.

\begin{algorithm}
\caption{CMTF --- Basic algorithm for the Tucker format}\label{agm:CMTF-Tucker}
\renewcommand{\algorithmicrequire}{\textbf{Input:}}
\renewcommand{\algorithmicensure}{\textbf{Output:}}
\begin{algorithmic}
\Require $\calX\in\R^{m\times n_2\times n_3}$, $Y\in\R^{m\times n}$, $k<\min\{n_2,n_2,n\}$
\Ensure $\widetilde{\calX}\in\R^{m\times n_2\times n_3}$, $U\in\R^{m\times k}$, $W\in\R^{n\times k}$
\State $[U,V,W]=\texttt{CMF}(X_{(1)},Y,k)$ \Comment{Algorithm~\ref{agm:CMF}}
\State $\widetilde{X}_{(1)}=UV^T$
\State Fold $\widetilde{X}_{(1)}$ into $\widetilde{\calX}$.
\end{algorithmic}
\end{algorithm}

In order to randomize Algorithm~\ref{agm:CMTF-Tucker}, one should use one of  Algorithms~\ref{agm:randCMF}, \ref{agm:rsiCMF}, or~\ref{agm:CMFRBKI} in place of  Algorithm~\ref{agm:CMF} (within  Algorithm~\ref{agm:CMTF-Tucker}).
Note that, compared to CMF, the sizes of the random Gaussian matrices created for the randomized CMTF will be larger. Instead of projection~\eqref{randproj}, we have
$$\widehat{X}_{(1)}=\Pi_{X_{(1)}\Omega}X_{(1)}, \quad \widehat{Y}=\Pi_{Y\Omega}Y.$$
That is, for $\calX\in\R^{m\times n_2\times n_3}$ and $Y\in\R^{m\times n}$, $\Omega_1$ is $n_2n_3\times k$ and $\Omega_2$ is $n\times k$ matrix.
Except for this difference,
the randomization process is the same.
Numerical results of the randomized algorithms are discussed in Section~\ref{sec:CMTF_num}.

We have only considered joint approximation coupled in the first mode. Instead, one can be interested in approximation coupled in a different way. In that case, relations~\eqref{CMTF_Tucker} and~\eqref{CMTF_TuckerU} are modified accordingly. For example, if we are looking for an approximation of the form
$$\calX\approx \calS\times_1V_1\times_2U\times_3V_3, \quad Y\approx UW^T,$$
we need to represent the tensor $\calX$ using mode-$2$ matricization,
$$X_{(2)}\approx U(\calS\times_1V_1\times_3V_3)_{(2)}.$$
This can be written as CMF of $X_{(2)}$ and $Y$, very similar to~\eqref{CMTF_Tucker2},
$$X_{(2)}\approx UV^T, \quad Y\approx UW^T,$$
which is then used in Algorithm~\ref{agm:CMTF-Tucker}.
For the coupling in the third mode, we proceed in the same way, but using the mode-$3$ matricization
$$X_{(3)}\approx U(\calS\times_1V_1\times_2V_2)_{(3)}$$
and CMF of $X_{(3)}$ and $Y$.
The coupling mode depends on the particular application and nature of the observed problem.

\subsection{Tensor in CP format}\label{ssec:CMTF_CP}

Next, we describe the CMTF problem where a tensor is represented in CP format. Again, let $\calX\in\R^{m\times n_2\times n_3}$ and $Y\in\R^{m\times n}$. Then,
\begin{equation}\label{CMTF_CP}
\calX\approx [[U,B,C]], \quad Y\approx UW^T,
\end{equation}
where $U\in\R^{m\times k}$, $B\in\R^{n_2\times k}$, $C\in\R^{n_3\times k}$, $W\in\R^{n\times k}$, is the joint rank-$k$ factorization of $\calX$ and $Y$ coupled in the first mode and taking $\calX$ in CP format, as it is shown in Figure~\ref{fig:CMTFCP}. Note that for CP format we consider tensor rank as it is given in the relation~\eqref{trank}.

\begin{figure}[t]
\begin{tikzpicture}
\draw [draw=black] (0,0) rectangle (2,2);
\node at (1,1) {$\calX$};
\draw (0,2)--(1,2.5);
\draw (2,2)--(3,2.5);
\draw (2,0)--(3,0.5);
\draw (1,2.5)--(3,2.5);
\draw (3,0.5)--(3,2.5);
\node at (3.5,1) {$\approx$};
\end{tikzpicture}
\begin{tikzpicture}
\draw [fill=gray, draw=gray] (0,0) rectangle (0.05,2);
\draw [draw=black] (0.1,1.95) rectangle (2.1,2);
\draw (0,2.05)--(0.05,2.05);
\draw (0.4,2.55)--(0.45,2.55);
\draw (0,2.05)--(0.4,2.55);
\draw (0.05,2.05)--(0.45,2.55);
\node at (2.6,1) {$+$};
\node at (0.35,0) {$u_1$};
\end{tikzpicture}
\begin{tikzpicture}
\draw [fill=gray, draw=gray] (0,0) rectangle (0.05,2);
\draw [draw=black] (0.1,1.95) rectangle (2.1,2);
\draw (0,2.05)--(0.05,2.05);
\draw (0.4,2.55)--(0.45,2.55);
\draw (0,2.05)--(0.4,2.55);
\draw (0.05,2.05)--(0.45,2.55);
\node at (2.9,1) {$+\cdots+$};
\node at (0.35,0) {$u_1$};
\end{tikzpicture}
\begin{tikzpicture}
\draw [fill=gray, draw=gray] (0,0) rectangle (0.05,2);
\draw [draw=black] (0.1,1.95) rectangle (2.1,2);
\draw (0,2.05)--(0.05,2.05);
\draw (0.4,2.55)--(0.45,2.55);
\draw (0,2.05)--(0.4,2.55);
\draw (0.05,2.05)--(0.45,2.55);
\node at (0.35,0) {$u_k$};
\end{tikzpicture}
\caption{Graphical description of the tensor $\calX$ from CMTF~\eqref{CMTF_CP}.}
\label{fig:CMTFCP}
\end{figure}

In order to find the coupled decomposition~\eqref{CMTF_CP} we first define the objective function
\begin{equation}\label{minCMTF}
f(U,B,C,W)=\|\calX-[[U,B,C]]\|^2+\|Y-UW^T\|_F^2\rightarrow\min.
\end{equation}
A common way for solving the optimization problem~\eqref{minCMTF} is the alternating least squares (ALS) method. That is an iterative method where one iteration is made of several, in this case four, microiterations. In each microiteration, optimization is done over one of the matrices $U,B,C,W$, while the others are fixed.
Let us see how the ALS algorithm for solving~\eqref{minCMTF} is developed.

We can write the function $f$ as
$$f(U,B,C,W)=f_1(U,B,C,W)+f_2(U,B,C,W),$$
where
\begin{align*}
f_1(U,B,C,W) & =\|\calX-[[U,B,C]]\|^2, \\
f_2(U,B,C,W) & =\|Y-UW^T\|_F^2. \end{align*}
It follows from~\cite[Corollary 4.2]{AcDuKo_chem11} that
\begin{align*}
\frac{\partial f_1}{\partial U} & = -2X_{(1)}(C\odot B) +2U\Gamma^{(1)}, \\
\frac{\partial f_1}{\partial B} & = -2X_{(2)}(C\odot U) +2B\Gamma^{(2)}, \\
\frac{\partial f_1}{\partial C} & = -2X_{(3)}(B\odot U) +2C\Gamma^{(3)}, \end{align*}
for
$$\Gamma^{(1)}=(B^TB)\ast (C^TC), \quad
\Gamma^{(2)}=(U^TU)\ast (C^TC), \quad
\Gamma^{(3)}=(U^TU)\ast (B^TB).$$
Obviously, $\frac{\partial f_1}{\partial W}=0$.
For the function $f_2$, after a short calculation, we get
$$\frac{\partial f_2}{\partial U}=-2YW+2UW^TW,\quad
\frac{\partial f_2}{\partial W}=-2Y^TU+2WU^TU,$$
while
$\frac{\partial f_2}{\partial B}=\frac{\partial f_2}{\partial C}=0$.
Furthermore, by setting the partial derivatives of $f$ to zero, we get
\begin{align*}
\frac{\partial f}{\partial U} & = -2X_{(1)}(C\odot B) +2U\Gamma^{(1)} -2YW+2UW^TW =0, \\
\frac{\partial f}{\partial B} & = -2X_{(2)}(C\odot U) +2B\Gamma^{(2)} =0, \\
\frac{\partial f}{\partial C} & = -2X_{(3)}(B\odot U) +2C\Gamma^{(3)} =0, \\
\frac{\partial f}{\partial W} & = -2Y^TU+2WU^TU =0.
\end{align*}
That is,
\begin{align*}
U & = (X_{(1)}(C\odot B)+YW)(\Gamma^{(1)}+W^TW)^{-1}, \\
B & = X_{(2)}(C\odot U)(\Gamma^{(2)})^{-1}, \\
C & = X_{(3)}(B\odot U)(\Gamma^{(3)})^{-1}, \\
W & = Y^TU(U^TU)^{-1}.
\end{align*}

Now we can write ALS algorithm for CMTF.

\begin{algorithm}
\caption{CMTF-CP ALS}\label{agm:CMTF-CPals}
\renewcommand{\algorithmicrequire}{\textbf{Input:}}
\renewcommand{\algorithmicensure}{\textbf{Output:}}
\begin{algorithmic}
\Require $X\in\R^{m\times n_2\times n_3}$, $Y\in\R^{m\times n}$, $k<\min\{n_2,n_3,n\}$
\Ensure $U\in\R^{m\times k}$, $B\in\R^{n_2\times k}$, $C\in\R^{n_3\times k}$, $W\in\R^{n\times k}$
\State Initialize $U,B,C,W$.
\Repeat
\State $U=X_{(1)}(C\odot B)+YW)(\Gamma^{(1)}+W^TW)^{-1}$
\State $B=X_{(2)}(C\odot U)(\Gamma^{(2)})^{-1}$
\State $C=X_{(3)}(B\odot U)(\Gamma^{(3)})^{-1}$
\State $W=Y^TU(U^TU)^{-1}$
\Until{convergence}
\end{algorithmic}
\end{algorithm}

Randomization of Algorithm~\ref{agm:CMTF-CPals} is more complex than for Algorithm~\ref{agm:CMF}. Let $\calX=[[A,B,C]]$ be a CP decomposition of $\calX\in\R^{m\times n_2\times n_3}$ and let $Y\in\R^{m\times n}$. First, we generate random Gaussian matrices  $\Omega_1\in{\R^{(n_2n_3)\times k}}$ and $\Omega_2\in{\R^{n\times k}}$. We compute thin QR decompositions $Q_1R_1=X_{(1)}\Omega_1$ and $Q_2R_2=Y\Omega_2$, such that $Q_1,Q_2\in\R^{m\times k}$, where $X_{(1)}$ is mode-$1$ matricization of $\calX$. Then, as it was explained in the Subsection~\ref{ssec:randCMF}, we reorthogonalize the columns of the matrix $\begin{bmatrix}
Q_1 & Q_2\end{bmatrix}$ to obtain $Q\in\R^{m\times (k+p)}$, where $p$ $(0\leq p\leq k)$ is an oversampling parameter. Randomized CMTF-CP ALS is presented in Algorithm~\ref{agm:randCMTF-CPals}.
Randomization using RSI and RBKI is done accordingly.

\begin{algorithm}
\caption{Randomized CMTF-CP ALS}\label{agm:randCMTF-CPals}
\renewcommand{\algorithmicrequire}{\textbf{Input:}}
\renewcommand{\algorithmicensure}{\textbf{Output:}}
\begin{algorithmic}
\Require $X\in\R^{m\times n_2\times n_3}$, $Y\in\R^{m\times n}$, $k<\min\{n_2,n_3,n\}$
\Ensure $U\in\R^{m\times k}$, $B\in\R^{n_2\times k}$, $C\in\R^{n_3\times k}$, $W\in\R^{n\times k}$
\State Generate random matrices $\Omega_1\in{\R^{(n_2\cdot n_3)\times k}}$ and $\Omega_2\in{\R^{n\times k}}$.
\State $[Q_1,\sim]=\texttt{qr}(X_{(1)}\Omega_1,0)$
\State $[Q_2,\sim]=\texttt{qr}(Y\Omega_2,0)$
\State $[Q,\sim]=\texttt{qr}([Q_1,Q_2],0)$
\State $[\widehat{U},B,C,W]=\texttt{CMTF}(\calX\times_1Q^T,Q^TY)$ \Comment{Algorithm~\ref{agm:CMTF-CPals}}
\State $U=Q\widehat{U}$
\end{algorithmic}
\end{algorithm}

Since Algorithm~\ref{agm:CMTF-CPals} is an iterative algorithm, convergence issues may occur. Therefore, it is worth observing other options.
Using the relation~\eqref{CPmat}, tensor $\calX$ from~\eqref{CMTF_CP} can be written as
$$X_{(1)}\approx U(C\odot B)^T.$$
Therefore, it makes sense to try to apply CMF on $X_{(1)}$ and $Y$, just as we do for the Tucker tensor format in~\eqref{CMTF_Tucker2}.
Hence, we can use Algorithm~\ref{agm:CMTF-Tucker} and we get
$$X_{(1)}\approx UV^T, \quad Y\approx UW^T.$$
Again, we cannot be sure that there are $B$ and $C$ of the appropriate sizes such that $V=C\odot B$. However, we get a joint factor $U$ and a numerical advantage of the approximation compared to the Algorithm~\ref{agm:CMTF-CPals} is seen in the next Section.

Similarly as for the Tucker tensor format, we have only examined coupling in the first mode. Results for the couplings in other modes follow the same reasoning, but use different matricizations of $\calX$. In the case of the second mode, we need
$$\calX\approx [[A,U,C]], \quad Y\approx UW^T,$$
and we use the relation
$$X_{(2)}=U(C\odot A)^T$$
to justify the use of CMF of $X_{(2)}$ and $Y$.
On the other hand, if we need
$$\calX\approx [[A,B,U]], \quad Y\approx UW^T,$$
the use of CMF of $X_{(3)}$ and $Y$ is a consequence of the relation
$$X_{(3)}=U(B\odot A)^T.$$

\section{Numerical examples for CMTF}\label{sec:CMTF_num}

Here we consider two numerical examples comparing our CMTF algorithms. In the first example, we are going to compare two approaches to CMTF algorithms without randomization, one using the ALS method and the other using SVD. Namely, we compare Algorithm~\ref{agm:CMTF-CPals} with Algorithm~\ref{agm:CMTF-Tucker}. The second example is analogous to Examples 1--5 from Section~\ref{sec:CMF_num} and we compare randomized versions of Algorithm~\ref{agm:CMTF-Tucker}.

\subsection{First example --- Comparing the basic algorithms}
Let $m=100$, $n_2=50$, $n_3=20$, $n=30$ and $r=3$ and define:
\begin{itemize}
\item \texttt{U = randn(m,r); B = randn(n2,r); C=randn(n3,r); W=randn(n,r);}
\item $\calX = [[U,B,C]]$;
\item  \texttt{Y=U*W';}
\end{itemize}

\begin{figure}[ht]
\begin{minipage}{.5\textwidth}
    \centering
    \includegraphics[width=1\linewidth]{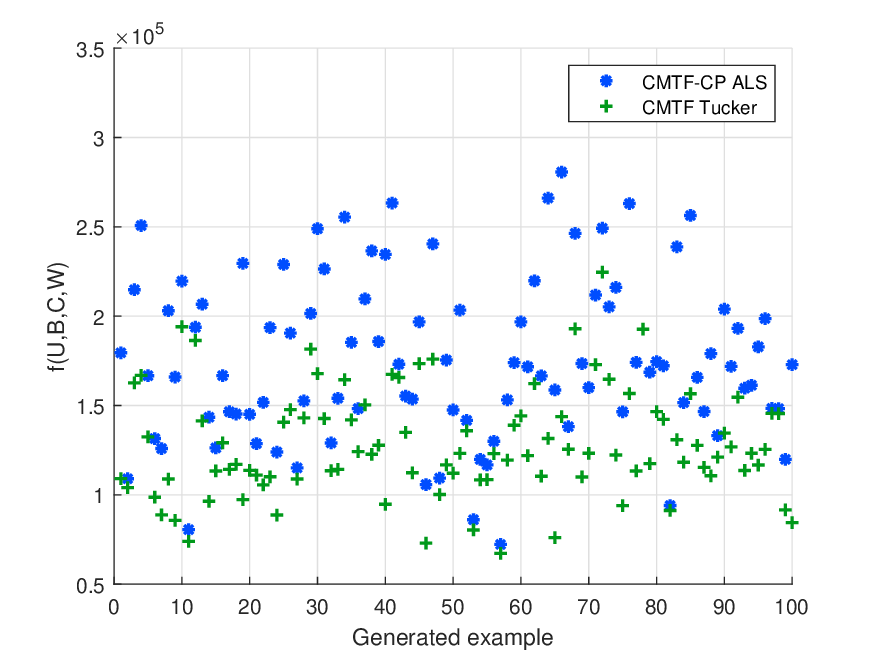}
    \end{minipage}%
     \begin{minipage}{.5\textwidth}
        \centering
        \includegraphics[width=1\textwidth]{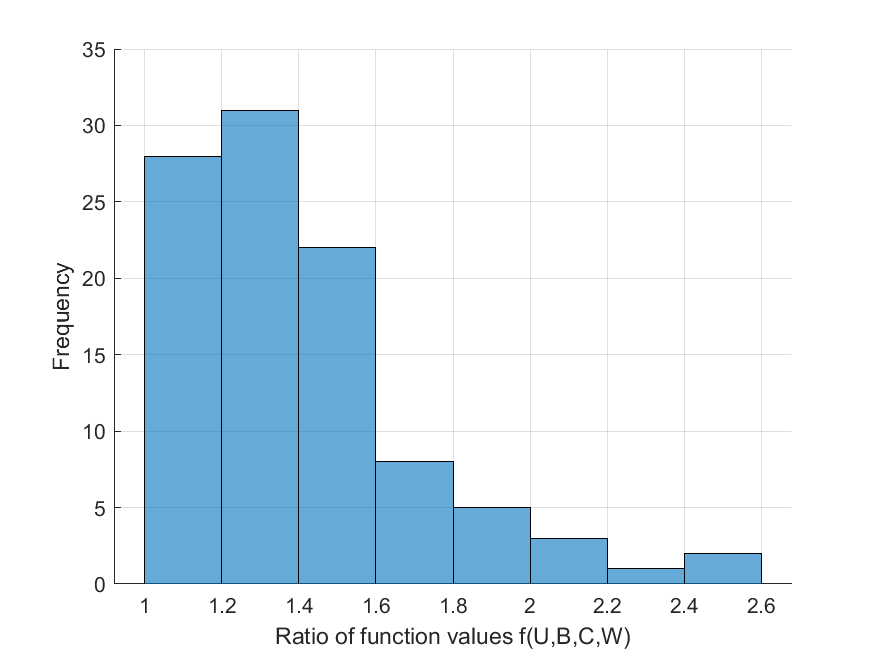}	
    \end{minipage}
    \caption{Comparison of Algorithms \ref{agm:CMTF-Tucker} and \ref{agm:CMTF-CPals} for CMTF for 100 randomly generated examples.}
    \label{fig:MainComparison}
\end{figure}

We use this example to test the algorithms without randomization, i.e., we compare Algorithm~\ref{agm:CMTF-Tucker} for a tensor in the Tucker format using CMF of $X_{(1)}$ and $Y$ and Algorithm~\ref{agm:CMTF-CPals} for a tensor in CP format using the ALS method. We generate $100$ random pairs $(\calX,Y)$ as described above. {The approximation rank is $k=3$.} For the values obtained by the algorithms, we compute
\begin{equation*}
    \|\mathcal{X}-\mathcal{S}\times_1 U \times_1 V_2 \times_3 V_3\|^2+\|Y-UW^T\|_F^2,
\end{equation*}
for $\calX$ in the Tucker format, and
\begin{equation*}
    \|\mathcal{X}-[[U,B,C]]\|^2+\|Y-UW^T\|_F^2,
\end{equation*}
for $\calX$ in CP format.
These expressions are the values of the corresponding objective functions.

The results are represented in Figure \ref{fig:MainComparison}. It is clear that Algorithm~\ref{agm:CMTF-Tucker}, that is, the method which uses CMF of the tensor matricization and matrix $Y$, performs equally well and almost always better than Algorithm~\ref{agm:CMTF-CPals}, which is based on ALS. In most experiments, the objective function value $f(U,B,C,W)$ obtained with ALS is at least $20\%$ higher, and in some cases, it is up to twice as high, compared to our approach. This is because convergence of the ALS method to the global minimum is not guaranteed, i.e., the method could converge to another stationary point. Therefore, in the next example, we do not use the ALS approach.

\subsection{Second example --- Synthetic tensor and matrix}
This example is inspired by the matrix $X$ from \texttt{TestSynthetic2}. We define a diagonal matrix
$$S = \operatorname{diag}(1,\ldots,1,d^{-2},d^{-3},\ldots,d^{-(n-r+1)})\in\R^{n\times n}$$
and construct matrix $Y\in\R^{n\times n}$ and tensor $X \in\R^{n\times n\times n}$ in the following way:
\begin{itemize}
    \item \texttt{UY=orth(rand(n));} \texttt{VY=orth(rand(n));}
    \item \texttt{Y = UY*S*VY';}
    \item \texttt{X(:,:,1)=[UY(:,1:r1) orth(rand(n,n-r1)]*S*[VY(:,1:r1) orth(rand(n,n-r1)]';}
    \item \texttt{X(:,:,2)=[UY(:,1:r2) orth(rand(n,n-r2)]*S*[VY(:,1:r2) orth(rand(n,n-r2)]';}
    \item \texttt{X(:,:,3)=[UY(:,1:r3) orth(rand(n,n-r3)]*S*[VY(:,1:r3) orth(rand(n,n-r3)]';}
\end{itemize}
For our experiment, we choose $n=100$, $r=5$, $d=2$, $r_1=5$, $r_2 = 10$ and $r_3=7$. We are looking for the rank $k=10$ approximation.
\begin{comment}
We separately analyze the algorithms based on the tensor format.
\end{comment}

\begin{table}[ht]
        \centering
        \begin{tabular}{|l|c|c|c|c|c|c|c|}
            \hline
            Algorithm  & $p$ & $\ell$ & $q$ & $err_{\calX}$ & $err_Y$ & \makecell{total\\ time (s)} & \makecell{CMF\\ time (s)} \\ \hline\hline
            Basic CMTF  & - & - & - & $2.51650855\cdot 10^{-2}$ & $2.52792977\cdot 10^{-2}$ & $0.444682$ & $0.444682$ \\
            Randomized  & 10 & - & - & $2.96659479\cdot 10^{-2}$ & $2.84378062\cdot 10^{-2}$ & $0.234335$ & $0.022847$ \\
            RSI  & 10 & - & 5 & $2.51760209\cdot 10^{-2}$ & $ 2.52882608\cdot 10^{-2} $ & $0.235490$ & $0.022847$ \\
            RBKI  & 26 & 1 & 18 & $2.51648934\cdot 10^{-2}$ & $2.52806919\cdot 10^{-2}$ & $0.247642$ & $0.012539$\\
            RBKI  & 22 & 2 & 8 & $2.51646244\cdot 10^{-2}$ & $ 2.52820672 \cdot 10^{-2}$ & $0.244876$ & $0.022847$\\
             \hline
        \end{tabular}
        \smallskip
        \caption{Selected results for CMTF where the tensor is in the Tucker format and $k=10$. Oversampling parameters, relative errors for RSI with $q=5$ and RBKI with $\ell=1,2$, and running times.}
        \label{tab:testdiffTensor}
\end{table}

\begin{figure}[ht]
\begin{multicols}{2}
    \includegraphics[width=0.9\linewidth]{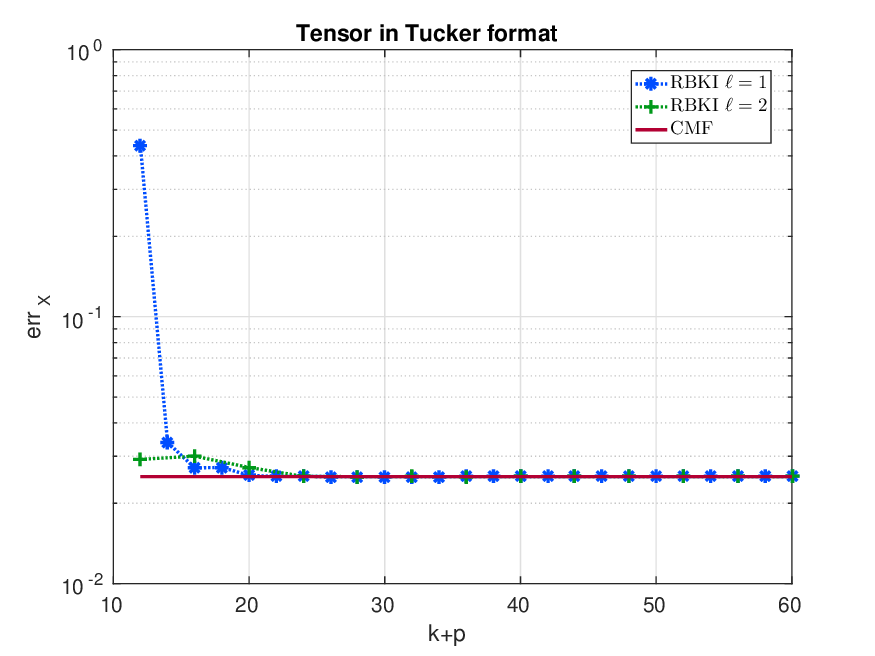}\par
    \includegraphics[width=0.9\linewidth]{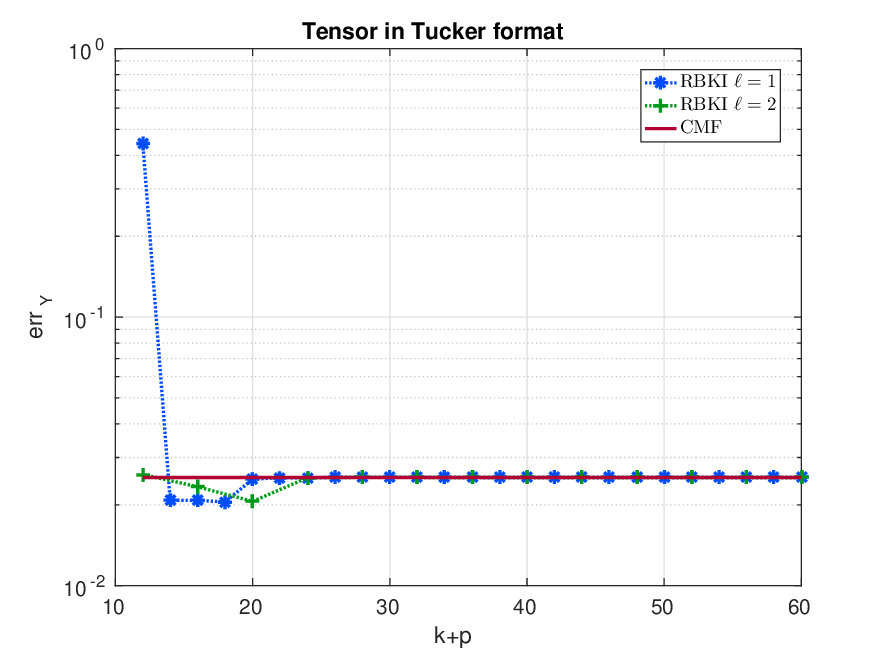}\par
    \end{multicols}
\caption{Relative errors with respect to number of iterations in RBKI algorithm for Matrix-Tensor decomposition in the Tucker format.}
\label{fig:RBKIMatrixTensor}
\end{figure}

Following the results from the previous example, here we do the tests only for the Tucker format. We use the result obtained by the basic Algorithm~\ref{agm:CMTF-Tucker} as a benchmark. In the randomized versions, we use Algorithm~\ref{agm:rsiCMF} for RSI and Algorithm~\ref{agm:RBKI} for RBKI, instead of  Algorithm~\ref{agm:CMF}.
The accuracy is determined by the relative errors,
\begin{equation}\label{eq:errCMTF-Tucker}
err_{\calX} = \frac{\|\calX - \calS\times_1U\times_2V_2\times_3V_3\|}{\|\calX\|}
\end{equation}
and $err_Y$ as it is given in~\eqref{eq:err}.
The selected results, together with the running times, are presented in Table~\ref{tab:testdiffTensor}. We can see that randomized algorithms are always faster, both in total time and CMF time.
We present Figure~\ref{fig:RBKIMatrixTensor} showing the relative errors for the RBKI method with $\ell=1$ and $\ell=2$, with various $q$.

\section{Face recognition}\label{sec:facerec}

So far we have introduced CMF and CMTF algorithms with different randomization techniques. Since these algorithms extract a common part of two matrices (or matrix and tensor), we created the new face recognition algorithms with CM(T)F as their main feature.

\begin{figure}[ht]
\centering
\includegraphics[width=\linewidth]{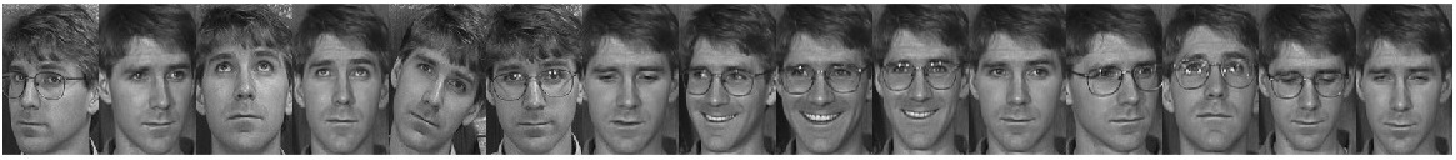}
\caption{An example of $15$ images per person from the Georgia Tech database. }\label{fig: FR_all_images}
\end{figure}

\begin{figure}[ht]
\centering
\includegraphics[width=0.65\linewidth]{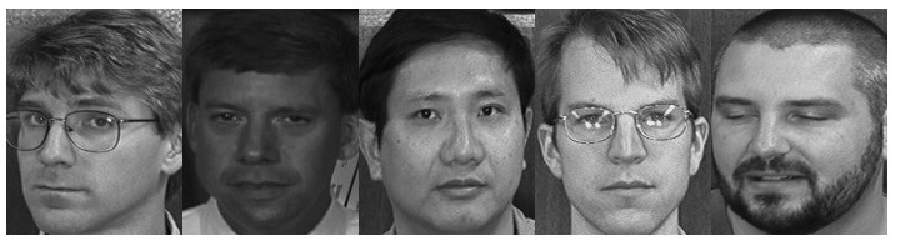}
\caption{Images of five people used in the tests of our face recognition algorithms.}\label{fig:5_persons}
\end{figure}

We performed the tests on the image collection from the widely used Georgia Tech database~\cite{GaTech}. This database contains a total of $750$ images of $50$ different people. Each person is represented by $15$ face images taken with different facial expressions.
An example of a set of images for one person is given in Figure~\ref{fig: FR_all_images}.
We used the first ten as the training images and the last five as the test images.
Developing a new optimal face recognition algorithm is out of the scope of this paper.
Therefore, to expedite testing, we selected $75$ images from five different individuals. Those are presented in Figure~\ref{fig:5_persons}. For now, we limit our tests to the grayscale images.
All images are uniformly cropped to the same size, with each image represented by a single matrix. Although the CMF algorithm requires only the same number of rows, the best test results were achieved when we also set the same number of columns.

We constructed two versions of the face recognition algorithm. The first one employs the coupled factorization of two matrices. The second one uses the coupled matrix and tensor factorization.
The approximation rank in CMF was set to $k=5$. Randomization parameters are $q=2$ for RSI and $\ell=5$, $q=2$ for RBKI. The algorithm steps are as follows.

\begin{enumerate}
\item[(1)] The first 10 images per person were put in our database, where it is known which person is in which image. We denoted the matrices representing those images with $X^{(1)}$, $X^{(2)}$, $\ldots,X^{(50)}$. Matrices $X^{(1)}$, $X^{(2)}$, $\ldots,X^{(10)}$ correspond to the first person, $X^{(11)}$, $X^{(12)}$, $\ldots,X^{(20)}$ to the second one, etc. We used the remaining 25 images (five per person) to test the algorithms. Matrix $Y$ represents the image of the person we want to recognize.

\item[(2 CMF)]
We calculate
\begin{equation*}
\operatorname{CMF}(X^{(1)},Y), \ \operatorname{CMF}(X^{(2)},Y), \ \ldots, \ \operatorname{CMF}( X^{(50)},Y).
\end{equation*}
The corresponding errors are calculated as
$$err^{(i)}=err_X^{(i)}+err_Y^{(i)}, \quad \text{for all } i=1,\ldots,50,$$
where $err_X$ and $err_Y$ are the relative errors defined in the relation~\eqref{eq:err}.

\item[(2 CMTF)]
We created five $m\times n\times10$ tensors, one per person, each of them containing ten images of the same person.
Then, we calculate
\begin{equation*}
\operatorname{CMTF}(\calX^{(1)},Y), \ \operatorname{CMTF}(\calX^{(2)},Y), \ \ldots, \ \operatorname{CMTF}(\calX^{(5)},Y).
\end{equation*}
The relative errors are determined by
$$err^{(i)}=err_{\calX}^{(i)}+err_Y^{(i)}, \quad \text{for all } i=1,\ldots,5,$$
where $err_{\calX}$ is like in the relation~\eqref{eq:errCMTF-Tucker}, for tensors in the Tucker format, or
$$err_{\calX} = \frac{\|\calX - [[U,B,C]]\|}{\|\calX\|},$$
in CP format, and $err_Y$ is as in the relation~\eqref{eq:err}.
\item[(3)] The index $i$ for which the value $err^{(i)}$ is minimal determines which person is in the image.
\end{enumerate}

We compared $12$ different versions of the algorithm, depending on whether CMF or CMTF is used and also on the randomization type (basic algorithm without randomization, simple randomization, RSI, or RBKI). For each tested case, we calculated the success rate, that is, the number of correctly recognized images divided by the total number of tested images. The results are presented in Table~\ref{tab:face_rec}.

\begin{table}[ht]
        \centering
        \begin{tabular}{|l|c|c|c|c|c||c|}
            \hline
            Algorithm &  Person 1 & Person 2 & Person 3 & Person 4 & Person 5 & Total  \\
             & (\%) & (\%) & (\%)& (\%)& (\%) & (\%) \\
            \hline\hline
            CMF  & $80 $  & $100 $  & $20 $ &  $100 $ &  $100 $  & $80 $ \\
            Randomized CMF   & $80 $  & $100 $  & $40 $ & $80 $ & $80 $ &  $76 $   \\
            RSI CMF  & $100 $  & $100 $  &  $20 $ &  $80 $ &  $100 $ &   $80 $  \\
            RBKI CMF  & $80 $  & $100 $  & $20 $ &  $100 $ &  $100 $  & $80 $ \\
            \hline
            CMTF-Tucker  & $60 $  & $100 $  & $100 $ &  $100 $ &  $80 $  & $88 $ \\
             Randomized CMTF-Tucker   & $40 $  & $80 $  & $20 $ & $80 $ & $20 $ &  $48 $   \\
            RSI CMTF-Tucker  & $60 $  & $100 $  &  $80 $ &  $100 $ &  $80 $ &   $84 $  \\
            RBKI CMTF-Tucker  & $60 $  & $100 $  & $100 $ &  $100 $ &  $80 $  & $88 $ \\
            \hline
             CMTF-CP ALS  & $80 $  & $100 $  & $0 $ &  $40 $ &  $20 $  & $48 $ \\
             Randomized CMTF-CP ALS   & $40 $  & $100 $  & $0 $ & $100 $ & $0 $ &  $48 $   \\
             RSI CMTF-CP ALS  & $0 $  & $100 $  &  $0 $ &  $100 $ &  $60 $ &   $52 $  \\
             RBKI CMTF-CP ALS  & $20 $  & $100 $  & $40 $ &  $100 $ &  $60 $  & $64 $ \\
            % \hline
            % CMTF-CP & $60 $  & $100 $  & $100 $ &  $100 $ &  $80 $  & $88 $ \\
            % Randomized CMTF-CP & $40 $  & $80 $  & $20 $ & $80 $ & $20 $ &  $48 $   \\
            % RSI CMTF-CP & $60 $  & $100 $  &  $80 $ &  $100 $ &  $80 $ &   $84 $  \\
            % RBKI CMTF-CP & $60 $  & $100 $  & $100 $ &  $100 $ &  $80 $  & $88 $ \\
            \hline
        \end{tabular}
        \smallskip
         \caption{The success rate for different versions of the face recognition algorithm.}
         \label{tab:face_rec}
    \end{table}

We can see that CMTF-Tucker shows the best performance. The success rate goes up to $88\%$. As expected, for the randomized algorithms, the best results are achieved when RBKI is used. For CMF and CMTF-Tucker, these are the same as for the basic Algorithms~\ref{agm:CMF} and~\ref{agm:CMTF-Tucker}.
With CMTF-Tucker, the success rate is more evenly distributed among individuals than with CMF.
The worst results were obtained with CMTF-CP ALS. Here, we can also observe that, unexpectedly,  randomization improved the basic algorithm. We attribute this to the previously mentioned convergence problems.

\section{Conclusion}\label{sec:conclusion}

In this paper, we analyze the coupled factorization of two matrices (CMF), as well as the coupled factorization of a matrix and a tensor (CMTF), and the possibilities for their randomization.
We showed that the problem of a rank-$k$ coupled matrix factorization for the matrices $X\in\R^{m\times n_1}$ and $Y\in\R^{m\times n_2}$ is equivalent to the problem of rank-$k$ approximation of the matrix $\begin{bmatrix}
X & Y
\end{bmatrix}\in\R^{m\times(n_1+n_2)}$. Moreover, we showed that a related problem of a rank-$k$ coupled factorization of tensor $\calX\in\R^{m\times n_2\times n_3}$ and matrix $Y\in\R^{m\times n}$ is equivalent to the problem of rank-$k$ approximation of the matrix $\begin{bmatrix}
X_{(1)} & Y
\end{bmatrix}\in\R^{m\times(n_2n_3+n)}$, where $X_{(1)}\in\R^{m\times(n_2n_3)}$ is mode-$1$ matricization of $\calX$.
The main part of the paper is concerned with randomization techniques that can be used with CMF and CMTF algorithms. First, we argue on the importance of choosing a good projection matrix for the dimension reduction of the original problem. We propose a method where we use QR decomposition of both terms in factorization, as opposed to just computing the QR decomposition of an augmented matrix. This produces more accurate results. Next, using our strategy for finding the projection matrix, and motivated by three randomization techniques --- randomization based on the randomized SVD algorithm, randomized subspace iteration (RSI), and randomized block Krylov iteration (RBKI), we created and tested different versions of randomized CM(T)F algorithms. Finally, we applied our algorithms to the face recognition problem, resulting in new competitive algorithms for this type of problem.

Matlab codes used to generate numerical examples are available at

\href{https://github.com/anitacarevic/Randomized-coupled-decompositions}{https://github.com/anitacarevic/Randomized-coupled-decompositions}.

\section*{Acknowledgments}
The authors are grateful to Daniel Kressner for his useful comments on the early version of this work. The authors also thank the anonymous referees for their helpful suggestions that improved the paper.

\bibliographystyle{siam}
\bibliography{ctd}

\begin{thebibliography}{10}

\bibitem{Acar_metabolimics15}
{\sc E.~Acar, R.~Bro, and A.~K. Smilde}, {\em Data fusion in metabolomics using coupled matrix and tensor factorizations}, Proceedings of the IEEE, 103 (2015), pp.~1602--1620.

\bibitem{AcDuKo_chem11}
{\sc E.~Acar, D.~M. Dunlavy, and T.~G. Kolda}, {\em A scalable optimization approach for fitting canonical tensor decompositions}, Journal of Chemometrics, 25 (2011), pp.~67--86.

\bibitem{AcDuKoMo11}
{\sc E.~Acar, D.~M. Dunlavy, T.~G. Kolda, and M.~Mørup}, {\em Scalable tensor factorizations for incomplete data}, Chemometrics and Intelligent Laboratory Systems, 106 (2011), pp.~41--56.
\newblock Multiway and Multiset Data Analysis.

\bibitem{Bagherian_bioinf20}
{\sc M.~Bagherian, R.~B. Kim, C.~Jiang, M.~A. Sartor, H.~Derksen, and K.~Najarian}, {\em {Coupled matrix–matrix and coupled tensor–matrix completion methods for predicting drug–target interactions}}, Briefings in Bioinformatics, 22 (2020), pp.~2161--2171.

\bibitem{Borsoi_bioinf23}
{\sc R.~A. Borsoi, I.~Lehmann, M.~A. B.~S. Akhonda, V.~D. Calhoun, K.~Usevich, D.~Brie, and T.~Adali}, {\em Coupled cp tensor decomposition with shared and distinct components for multi-task fmri data fusion}, in ICASSP 2023 - 2023 IEEE International Conference on Acoustics, Speech and Signal Processing (ICASSP), 2023, pp.~1--5.

\bibitem{DriMa16}
{\sc P.~Drineas and M.~W. Mahoney}, {\em Randnla: randomized numerical linear algebra}, Commun. ACM, 59 (2016), p.~80–90.

\bibitem{el2024tensor}
{\sc A.~El~Hachimi, K.~Jbilou, A.~Ratnani, and L.~Reichel}, {\em A tensor bidiagonalization method for higher-order singular value decomposition with applications}, Numerical Linear Algebra with Applications, 31 (2024), p.~e2530.

\bibitem{hached2021multidimensional}
{\sc M.~Hached, K.~Jbilou, C.~Koukouvinos, and M.~Mitrouli}, {\em A multidimensional principal component analysis via the c-product golub--kahan--svd for classification and face recognition}, Mathematics, 9 (2021), p.~1249.

\bibitem{HaMaTr11}
{\sc N.~Halko, P.~G. Martinsson, and J.~A. Tropp}, {\em Finding structure with randomness: probabilistic algorithms for constructing approximate matrix decompositions}, SIAM Rev., 53 (2011), pp.~217--288.

\bibitem{jia2017color}
{\sc Z.-G. Jia, S.-T. Ling, and M.-X. Zhao}, {\em Color two-dimensional principal component analysis for face recognition based on quaternion model}, in Intelligent Computing Theories and Application: 13th International Conference, ICIC 2017, Liverpool, UK, August 7-10, 2017, Proceedings, Part I 13, Springer, 2017, pp.~177--189.

\bibitem{KiersCP}
{\sc H.~A.~L. Kiers}, {\em Towards a standardized notation and terminology in multiway analysis}, Journal of Chemometrics, 14 (2000), pp.~105--122.

\bibitem{KoBa09}
{\sc T.~G. Kolda and B.~W. Bader}, {\em Tensor decompositions and applications}, SIAM Rev., 51 (2009), pp.~455--500.

\bibitem{Acar_metabolimics24b}
{\sc L.~Li, S.~Yan, D.~Horner, M.~Rasmussen, A.~Smilde, and E.~{Acar Ataman}}, {\em Revealing static and dynamic biomarkers from postprandial metabolomics data through coupled matrix and tensor factorizations}, Metabolomics, 20 (2024).

\bibitem{Ma16}
{\sc P.-G. Martinsson}, {\em Randomized methods for matrix computations}, IAS/Park City Mathematics Series,  (2016).

\bibitem{MaTr20}
{\sc P.-G. Martinsson and J.~A. Tropp}, {\em Randomized numerical linear algebra: Foundations and algorithms}, Acta Numerica, 29 (2020), p.~403–572.

\bibitem{meyer2024unreasonable}
{\sc R.~Meyer, C.~Musco, and C.~Musco}, {\em On the unreasonable effectiveness of single vector krylov methods for low-rank approximation}, in Proceedings of the 2024 Annual ACM-SIAM Symposium on Discrete Algorithms (SODA), SIAM, 2024, pp.~811--845.

\bibitem{MurrayRandNLA23}
{\sc R.~Murray, J.~Demmel, M.~W. Mahoney, N.~B. Erichson, M.~Melnichenko, O.~A. Malik, L.~Grigori, P.~Luszczek, M.~Derezinski, M.~E. Lopes, T.~Liang, H.~Luo, and J.~Dongarra}, {\em Randomized numerical linear algebra : A perspective on the field with an eye to software}, arXiv:2302.11474 [math.NA],  (2023).

\bibitem{nakatsukasa2020fast}
{\sc Y.~Nakatsukasa}, {\em Fast and stable randomized low-rank matrix approximation}, arXiv:2009.11392 [math.NA],  (2020).

\bibitem{GaTech}
{\sc A.~V. Nefian}, {\em Georgia {T}ech face database}.
\newblock \url{http://www.anefian.com/research/face_reco.htm}.

\bibitem{saibaba2023randomized}
{\sc A.~K. Saibaba and A.~Mi{{e}}dlar}, {\em Randomized low-rank approximations beyond gaussian random matrices}, arXiv:2308.05814 [math.NA],  (2023).

\bibitem{Schenker21}
{\sc C.~Schenker, J.~E. Cohen, and E.~Acar}, {\em A flexible optimization framework for regularized matrix-tensor factorizations with linear couplings}, IEEE Journal of Selected Topics in Signal Processing, 15 (2021), pp.~506--521.

\bibitem{Schenker21proc}
\leavevmode\vrule height 2pt depth -1.6pt width 23pt, {\em An optimization framework for regularized linearly coupled matrix-tensor factorization}, in 2020 28th European Signal Processing Conference (EUSIPCO), 2021, pp.~985--989.

\bibitem{SiGoCMF08}
{\sc A.~P. Singh and G.~J. Gordon}, {\em Relational learning via collective matrix factorization}, in Proceedings of the 14th ACM SIGKDD International Conference on Knowledge Discovery and Data Mining, KDD '08, New York, NY, USA, 2008, Association for Computing Machinery, p.~650–658.

\bibitem{Smilde_chem00}
{\sc A.~K. Smilde, J.~A. Westerhuis, and R.~Boqué}, {\em Multiway multiblock component and covariates regression models}, Journal of Chemometrics, 14 (2000), pp.~301--331.

\bibitem{SoDeL15-2}
{\sc M.~S\o{r}ensen, I.~Domanov, and L.~De~Lathauwer}, {\em Coupled canonical polyadic decompositions and (coupled) decompositions in multilinear rank-{$(L_{r,n},L_{r,n},1)$} terms---{P}art {II}: {A}lgorithms}, SIAM J. Matrix Anal. Appl., 36 (2015), pp.~1015--1045.

\bibitem{Armstrong_chrom23}
{\sc M.~D. {Sorochan Armstrong}, J.~L. Hinrich, A.~P. {de la Mata}, and J.~J. Harynuk}, {\em Parafac2×n: Coupled decomposition of multi-modal data with drift in n modes}, Analytica Chimica Acta, 1249 (2023), p.~340909.

\bibitem{SoDeL13proc}
{\sc M.~Sørensen and L.~De~Lathauwer}, {\em Coupled tensor decompositions for applications in array signal processing}, in 2013 5th IEEE International Workshop on Computational Advances in Multi-Sensor Adaptive Processing (CAMSAP), 2013, pp.~228--231.

\bibitem{TrWe23}
{\sc J.~A. Tropp and R.~J. Webber}, {\em Randomized algorithms for low-rank matrix approximation: Design, analysis, and applications}, arXiv:2306.12418 [math.NA],  (2023).

\bibitem{Tucker66}
{\sc L.~R. Tucker}, {\em Some mathematical notes on three-mode factor analysis}, Psychometrika, 31 (1966), pp.~279--311.

\bibitem{EeDeLChem20}
{\sc F.~{Van Eeghem} and L.~{De Lathauwer}}, {\em Tensor similarity in chemometrics}, in Comprehensive Chemometrics (Second Edition), S.~Brown, R.~Tauler, and B.~Walczak, eds., Oxford, 2020, Elsevier, pp.~337--354.

\bibitem{Smilde10}
{\sc I.~{Van Mechelen} and A.~K. Smilde}, {\em A generic linked-mode decomposition model for data fusion}, Chemometrics and Intelligent Laboratory Systems, 104 (2010), pp.~83--94.
\newblock OMICS.

\bibitem{vasilescu2002multilinear}
{\sc M.~A.~O. Vasilescu and D.~Terzopoulos}, {\em Multilinear analysis of image ensembles: Tensorfaces}, in Computer Vision—ECCV 2002: 7th European Conference on Computer Vision Copenhagen, Denmark, May 28--31, 2002 Proceedings, Part I 7, Springer, 2002, pp.~447--460.

\bibitem{wei2016tikhonov}
{\sc Y.~Wei, P.~Xie, and L.~Zhang}, {\em Tikhonov regularization and randomized gsvd}, SIAM Journal on Matrix Analysis and Applications, 37 (2016), pp.~649--675.

\bibitem{Woo14}
{\sc D.~P. Woodruff}, {\em Sketching as a tool for numerical linear algebra}, Found. Trends Theor. Comput. Sci., 10 (2014), p.~1–157.

\bibitem{xiao2018two}
{\sc X.~Xiao and Y.~Zhou}, {\em Two-dimensional quaternion pca and sparse pca}, IEEE transactions on neural networks and learning systems, 30 (2018), pp.~2028--2042.

\bibitem{Acar_metabolimics24a}
{\sc S.~Yan, L.~Li, D.~Horner, P.~Ebrahimi, B.~Chawes, L.~Dragsted, M.~Rasmussen, A.~Smilde, and E.~{Acar Ataman}}, {\em Characterizing human postprandial metabolic response using multiway data analysis}, Metabolomics, 20 (2024).

\bibitem{yang2004two}
{\sc J.~Yang, D.~Zhang, A.~F. Frangi, and J.-y. Yang}, {\em Two-dimensional pca: a new approach to appearance-based face representation and recognition}, IEEE transactions on pattern analysis and machine intelligence, 26 (2004), pp.~131--137.

\end{thebibliography}

\end{document}